\documentclass[12pt]{article}

\usepackage{amssymb}
\usepackage{amsmath}
\usepackage{graphicx}
\usepackage{color}

\newtheorem{thm}{Theorem}[section]
 \newtheorem{cor}[thm]{Corollary}
 \newtheorem{lem}[thm]{Lemma}
 \newtheorem{prop}[thm]{Proposition}
 \newtheorem{defn}[thm]{Definition}
 \newtheorem{rem}[thm]{Remark}
 \newtheorem{ex}[thm]{Example}
 \numberwithin{equation}{section}

\usepackage[shortlabels]{enumitem}
\newcommand{\qed}{\hfill{\rule{1ex}{1ex}}}
\DeclareMathOperator{\ind}{ind}

\DeclareMathOperator{\diag}{diag}
\DeclareMathOperator{\wind}{wind}

\DeclareMathOperator{\im}{im}
\DeclareMathOperator{\coker}{coker}

\usepackage{amssymb}
\usepackage{amsmath}
\usepackage{color}

\usepackage{enumitem}

 \newcommand{\vp}{\varphi}
 \newcommand{\ve}{\varepsilon}

\renewcommand{\Im}{\mathsf{Im}\,}

\newcommand{\nn}{\nonumber}

\newcommand{\bA}{{\bf A}}
\newcommand{\bB}{{\bf B}}

\newcommand{\bD}{{\bf D}}

\newcommand{\cA}{\mathcal{A}}
\newcommand{\cB}{\mathcal{B}}

\newcommand{\cF}{\mathcal{F}}

\newcommand{\cL}{\mathcal{L}}

\newcommand{\cP}{\mathcal{P}}
\newcommand{\cQ}{\mathcal{Q}}

\newcommand{\fN}{{\mathfrak N}}

\newcommand{\sC}{{\mathbb C}}
\newcommand{\sD}{{\mathbb D}}

\newcommand{\sN}{{\mathbb N}}

\newcommand{\sR}{{\mathbb R}}
\newcommand{\sS}{{\mathbb S}}
\newcommand{\sT}{{\mathbb T}}

\newcommand{\sW}{{\mathbb W}}

\newcommand{\sZ}{{\mathbb Z}}

\begin{document}

\vspace*{10mm}

\begin{center}
{\Large\textbf{Invertibility Issues for Toeplitz plus
Hankel\\[0.5ex] Operators and Their Close Relatives}}\footnote{This work was  supported by
the Special Project on High-Performance
Computing of the National Key R\&D Program of China (Grant No.~2016YFB0200604), the National Natural Science Foundation of China (Grant No.~11731006) and the Science Challenge Project of China (Grant No.~TZ2018001)..}
\end{center}

\vspace{5mm}

\begin{center}

\textbf{Victor D. Didenko and Bernd Silbermann}

\end{center}

\vspace{2mm}

 \begin{flushleft}
\small{SUSTech International Center for Mathematics, Southern University of Science and Technology, Shenzhen, China; \texttt{diviol@gmail.com}

\vspace{2mm}

Technische Universit{\"a}t Chemnitz, Fakult{\"a}t f\"ur Mathematik,
09107 Chemnitz, Germany; \texttt{silbermn@mathematik.tu-chemnitz.de}}

 \end{flushleft}

  \vspace{5mm}

\noindent \textbf{2010 Mathematics Subject Classification:} Primary 47B35, 47B38; Secondary 47B33,  45E10

\vspace{1mm}

\noindent\textbf{Key Words:} Toeplitz plus Hankel operators, Wiener-Hopf plus Hankel operators, invertibility, inverses

\vspace{10mm}

 \setcounter{page}{1} \thispagestyle{empty}

\begin{abstract}
The paper describes various approaches to the invertibility of Toe\-plitz plus Hankel operators in Hardy and  $l^p$-spaces, integral and difference Wiener-Hopf plus Hankel operators and generalized Toeplitz plus Hankel operators. Special attention is paid to a newly developed method, which allows to establish necessary, sufficient and also necessary and sufficient conditions of invertibility, one-sided and generalized invertibility for wide classes of operators and derive efficient formulas for the corresponding inverses. The work also contains a number of problems whose solution would be of interest in both theoretical and applied contexts.
\end{abstract}

%\tableofcontents % for long articles

%%%%%%%%%%%%%%%%%%%%%%%%%%%%%%%%%%%%%%%%%%%%%

%%% ----------------------------------------------------------------------
\tableofcontents

 \section{Introduction.\label{sec0}}
Toeplitz $T(a)$ and Hankel $H(b)$ operators appear in various fields of mathematics, physics, statistical mechanics and they have been thoroughly studied \cite{BS:2006,Kr:1987,Pe:2003}. Toeplitz plus Hankel operators $T(a)+H(b)$ and  Wiener-Hopf plus Hankel operators $W(a)+H(b)$ play an important role in random matrix theory  \cite{BR:2001,BCE:2010,FF:2004} and scattering  theory \cite{GR:2014, GR:2015, GR:2018,MST:1992,LMT:1992,MPST:1992,Te:1991}. Although Fredholm  properties  and index formulas for such operators acting on different Banach and Hilbert spaces are often known  --- cf.~\cite{BS:2006,DS:2013,KS:2001, KS:2000,KLR:2007,KLR:2009,Po:1979,RSS:2011,RS:1990,RS:2012, Si:1987}, their invertibility is little studied. So far progress has been made only in rare special cases. In this work, we want to present an approach, which allows to treat invertibility problem for a wide classes of Toeplitz plus Hankel operators on classic Hardy spaces and also for their close relatives: Toeplitz plus Hankel operators on $l^p$-spaces, generalized Toeplitz plus Hankel operators and integral and difference Wiener-Hopf plus Hankel operators.  The operators acting on classical Hardy spaces are discussed in more details, whereas for other classes of operators we only provide a brief overview of the corresponding results.

\section{Toeplitz and Hankel Operators on Hardy Spaces.\label{sec2}}

Let  $\sT:=\{t\in \sC: |t|=1\}$ be the counterclockwise oriented unit
circle in the complex plane $\sC$ and let $p\in [1,\infty]$. Consider the Hardy spaces
\begin{align*}
 H^p&=H^p(\sT):= \{f\in L^p(\sT): \widehat{f}_n=0\quad\text{for all}\quad n<0\},\\
\overline{H^p}&=\overline{H^p}(\sT):=\{f\in L^p(\sT): \widehat{f}_n=0\quad\text{for all}\quad n>0\},
\end{align*}
where $\widehat{f}_n$, $n\in \sN$ are the Fourier coefficients of  function $f$. Moreover, let $I$ denote the identity operator, $P: L^p(\sT)\to H^p(\sT)$ the  projection defined by
\begin{equation*}%\label{Eq}
P: \sum_{n=-\infty}^{\infty} \widehat{f}_n e^{in\theta} \mapsto
\sum_{n=0}^{\infty} \widehat{f}_n e^{in\theta}
\end{equation*}
and $Q=I-P$. If $p\in (1,\infty)$, the Riesz projection $P$ is bounded and $\im P= H^p$. 

On the space $H^p$, $1<p<\infty$, any function $a\in L^\infty$ generates two operators --- viz. the Toeplitz operator $T(a)\colon :f\mapsto Paf$ and the
Hankel operator $H(a)\colon f\mapsto PaQJf$, where $J:L^p\mapsto L^p$ is the  flip operator,
 \begin{equation*}
(Jf)(t):=t^{-1}f\left (t^{-1}\right ), \quad t\in \sT.
 \end{equation*}
We note that
 $$
 J^2=I,\; JPJ=Q,\; JQJ=P,\; JaJ=\widetilde{a},\; \widetilde{a}(t):=a(1/t),
 $$
and  the operators $T(a)$ and $H(b)$ are related to each other as follows
\begin{equation}\label{Widom}
\begin{aligned}
 T(ab) & = T(a)T(b)+H(a)H(\widetilde{b}),  \\
   H(ab)&=T(a)H(b)+H(a)T(\widetilde{b}).
\end{aligned}
\end{equation}

We now consider the invertibility of  Toeplitz plus Hankel operators $T(a)+H(b)$ generated by $L^\infty$-functions $a$ and $b$ and acting on a Hardy space $H^p(\sT)$. Observe that the matrix representation of such operators in the standard basis $\{t^n\}_{n=0}^\infty$ of $H^p(\sT)$ is
\begin{equation*}
T(a)+H(b) \sim(\widehat{a}_{k-j}+\widehat{b}_{k+j+1})_{k,j=0}^\infty.
\end{equation*}
There are a variety of approaches to the study of their invertibility and we briefly discuss some of them.

\subsection{Classical  approach: I. Gohberg, N. Krupnik and G. Litvinchuk.}

Let $\mathcal{L}(X)$ and $\cF(X)$ be, respectively, the sets of  linear bounded and Fredholm operators on the Banach space $X$. Besides, if $\cA$ is a unital Banach algebra, then $G\cA$ stands for the group of all invertible elements in $\cA$.

Assume that $a \in G L^\infty$, $b \in L^\infty$ and set
\begin{equation*}
V(a,b):=
\left(%
\begin{array}{cc}
 a- b\widetilde{b}\widetilde{a}^{-1} & d \\
 -c   &  \widetilde{a}^{-1}\\
   \end{array}%
\right),
\end{equation*}
where  $c:=\widetilde{b}\widetilde{a}^{-1}$,
$d:=b\widetilde{a}^{-1}$. Writing $T(a)\pm H(b)+Q$ for $(T(a)\pm H(b))P+Q$, we consider the operator $\diag(T(a)+H(b)+Q,T(a)-H(b)+Q)\colon L^p(\sT)\times L^p(\sT)$ and represent it in the form
 \begin{align}
 & \left(%
\begin{array}{cc}
  T(a)+H(b)+Q &  0\\
   0 &  T(a)-H(b)+Q\\
   \end{array}%
\right)\nn\\
&\hspace*{42mm}=
 A(T(V(a,b))+\diag (Q,Q))B\label{eqTV}
\end{align}
with invertible operators $A$ and $B$. More precisely,
\begin{equation*}
B= \left(%
\begin{array}{cc}
  I &  0\\
  \widetilde{b}I  &  \widetilde{a}I\\
   \end{array}%
\right)
\left(%
\begin{array}{cc}
  I & I\\
  J  &  -J\\
   \end{array}%
\right).
\end{equation*}
and the operator $A$ is also known but its concrete form is not important now.

An immediate consequence of the Eq.~\eqref{eqTV} is that both operators $T(a)\pm H(b)$ are Fredholm if and only if the block Toeplitz operator $T(V(a,b))$ is Fredholm. This representation is  of restricted use because  there are piecewise continuous functions $a,b$ such that only one of the operators $T(a)\pm H(b)$ is Fredholm. In addition, both operators $T(a)\pm H(b)$ can be Fredholm but may have different indices. Therefore, more efficient methods for studying Toeplitz plus Hankel operators are needed, especially for discontinuous generating functions $a$ and $b$.

Let us recall that there is a well-developed Fredholm theory for the operators $T(a)+H(b)$ with generating functions $a,b$ from the set of piecewise continuous functions $PC= PC(\sT)$ --- cf.~\cite{BS:2006} for operators acting on the space $H^2$ and \cite{RS:2012} for the ones on $H^p$, $p\neq 2$. However, the defect numbers of these operators, conditions for their invertibility, and inverse operators can be rarely determined directly from the Eq. \eqref{eqTV}.

\subsection{Basor-Ehrhardt approach.}

This approach is aimed at the study of defect numbers of $T(a)+H(b)\in \cF(H^p)$ by means of a factorization theory. It is well-known that for $b=0$, the problem can be solved by Wiener-Hopf factorization. Since this notion is important for what follows, we recall the definition here. Note that from now on, all operators are considered in the spaces $L^p$ or $H^p$ for $p\in (1,\infty)$. Moreover, let us agree that whenever $p$ and $q$ appear in the text, they are related as $1/p+1/q=1$.
\begin{defn}
We say that a function $a\in L^\infty$ admits a  Wiener-Hopf
factorization in $L^p$ if it can be represented in the form
 \begin{equation}\label{fac}
    a=a_- \chi_m a_+,
\end{equation}
where $a_+\in H^q$, $a_+^{-1}\in H^p$, $a_-\in \overline{H^p}$, $a_-^{-1}\in \overline{H^q}$, $\chi_m(t):=t^m$, $m\in \sZ$,
the term  $a_+^{-1}P(a_+ \vp)$ belongs to  $L^p$ for any $\vp$ from the set of all trigonometrical polynomials $\cP=\cP(\sT)$ and there is a constant $c_p$ such that
\begin{equation*}%\label{eqn1b}
||a_+^{-1}P(a_+ \vp)||_p\leq c_p ||\vp||_p \quad \text{for all} \quad \vp\in \cP.
\end{equation*}
 \end{defn}

 \begin{thm}[Simonenko \cite{Sim:1968}]
  If $a\in L^\infty$, then $T(a)\in \cF(H^p)$ if and only if $a\in GL^\infty$ and admits the Wiener-Hopf factorization \eqref{fac} in $L^p$. In this case
   \begin{equation*}
\ind T(a)=-m.
   \end{equation*}
   \end{thm}
 This result extends to the case of matrix-valued generating functions as follows.
  \begin{thm}
  If $a\in L^\infty_{N\times N}$, then $T(a)\in \cF(H^p_N)$ if and only if $a\in GL^\infty_{N\times N}$ and admits a factorization $a=a_- d a_+$, where
  \begin{align*}
&a_+\in H^q_{N\times N}, \; a_+^{-1}\in H^p_{N\times N},\; a_-\in
\overline{H}^p_{N\times N}, \; a_-^{-1}\in \overline{H}^q_{N\times N},\\
&d=\diag (\chi_{k_1},\chi_{k_2},\cdots,\chi_{k_N}), \quad \kappa_1,\kappa_2,\cdots, \kappa_N\in \sZ,
  \end{align*}
 the term  $a_+^{-1}P(a_+ \vp)$ belongs to  $L^p_N$ for any $\vp\in \cP_N$ and there is a constant $c_p$ such that
\begin{equation*}
||a_+^{-1}P(a_+ \vp)||_p\leq c_p ||\vp||_p \quad \text{for all} \quad \vp\in \cP_N.
\end{equation*}
Moreover, if $T(a)\in \cF(H^p_N)$, then
\begin{equation*}
\dim\ker T(a)=-\sum_{\kappa_j<0} \kappa_j, \quad \dim\coker T(a)=-\sum_{\kappa_j>0} \kappa_j.
\end{equation*}
   \end{thm}

The numbers $\kappa_j$, called partial indices, are uniquely defined. Moreover, in some sense, the Wiener-Hopf factorization is unique if it exists. For example, for $N=1$ the uniqueness of factorization can be ensured by the condition $a_-(\infty)=1$. We also note that if $T(a)\in \cF(H^p)$ and $\ind T(a)\geq 0$ ($\ind T(a)\leq 0$) then $T(a)$ is right-invertible (left-invertible) and if  $\kappa:=\ind T(a)>0$, the functions
\begin{equation*}
 a_+^{-1}\chi_j, \quad j=0,\ldots, \kappa-1
\end{equation*}
form a basis in the space $\ker T(a)$ and one of the right-inverses has the form $T^{-1}(a\chi_\kappa)T(\chi_\kappa)$.

 A comprehensive information about Wiener-Hopf factorization is provided in \cite{CG:1981, LS:1987} and in books, which deal with singular integral and convolution operators \cite{BS:2006, GK:1992a,GK:1992a, GK:1992b}. In particular, Wiener-Hopf factorization furnishes conditions for invertibility of related operators. However, generally there are no efficient methods for constructing such factorizations and computing partial indices even for continuous matrix-functions. Therefore, in order to study the invertibility of Toeplitz plus Hankel operators, we have to restrict ourselves to suitable classes of generating functions.

 In the beginning of this century, Ehrhardt \cite{E:2004h,E:2004} developed a factorization theory to study invertibility for large classes of convolution operators with flips. Toeplitz plus Hankel operators are included in this general framework. In particular, it was shown that an operator $T(a)+H(b)$, $a\in GL^\infty$, $b\in L^\infty$ is Fredholm if and only if the matrix-function
 \begin{equation*}
\phantom{aaaaa}V^\tau (a,b)=
\left(%
\begin{array}{cc}
 b \widetilde{a}^{-1} & a- b\widetilde{b}\widetilde{a}^{-1} \\
 \widetilde{a}^{-1}   &  -\widetilde{a}^{-1}\widetilde{b}\\
   \end{array}%
\right)
\end{equation*}
admits a certain type of antisymmetric factorization. Moreover, the defect numbers of the operator $T(a)+H(b)$ can be expressed via partial indices of this factorization. However, it is not known how the partial indices of general matrix-functions can be determined.

It was already mentioned that there are functions $a,b\in L^\infty$ such that the operator $T(a)+H(b)$ is Fredholm but $T(a)-H(b)$ is not. In this case we can use the representation \eqref{eqTV} and conclude that the matrix $V^\tau (a,b)$ does not admit a Wiener-Hopf factorization. Thus, if $V(a,b)$ admits a Wiener-Hopf factorization, then $V^\tau (a,b)$ has the antisymmetric factorization mentioned but not vice versa. %converse conclusion is not true.

This discussion  shows that one has to select a class of generating functions $a$ and $b$ such that the defect numbers of  operators $T(a)+H(b)\in \cF(H^p)$  can be determined. An important class of suitable generating functions $a\in GL^\infty(\sT)$, $b\in L^\infty(\sT)$ is given by the condition
\begin{equation}\label{eqn1c}
a \widetilde{a}=b \widetilde{b}.
\end{equation}
This class of pairs of functions first appears in \cite{BE:2013} and \cite{DS:2013}. The Eq.~\eqref{eqn1c} was latter called the matching condition and the corresponding duo $(a,b)$ a matching pair --- cf.~\cite{DS:2014}. Let us note that   Toeplitz plus Hankel operators of the form
\begin{equation}\label{eqn1d}
T(a)\pm H(a),\; T(a)- H(at^{-1}),\;T(a)+ H(at)
 \end{equation}
appear in random ensembles \cite{BR:2001,BCE:2010,FF:2004} and in numerical methods for singular integral equations on intervals \cite{JK:2019}. The generating functions of the operators \eqref{eqn1d} clearly satisfy the  matching condition \eqref{eqn1c}.

It is notable that the Fredholmness of the operator $T(a)+H(b)$ implies that $a\in GL^\infty(\sT)$. Therefore, the term  $b$ in the matching pair $(a,b)$ is also invertible in $L^\infty(\sT)$ and  one can introduce another pair $(c,d)$, called the subordinated pair for $(a,b)$, with the functions $c$ and $d$ defined by
\begin{equation*}
c:=a/b=\widetilde{b}/\widetilde{a},\; d:=a/\widetilde{b}=b/\widetilde{a}.
\end{equation*}
An important property of these functions is that $c \widetilde{c}=d \widetilde{d}=1$. In what follows,  any function $g\in L^\infty(\sT)$ satisfying the equation $g \widetilde{g}=1$  is referred to as a matching function.

Let us point out that the sets of matching functions and matching pairs are quite large. In particular, we have:
\begin{enumerate}
   \item Let $\sT^+:=\{t\in \sT:\Im t> 0 \}$ be
the upper half-circle. If an element $g_0\in GL^\infty$, then
 \begin{equation*}
g(t):=\left \{
 \begin{array}{ll}
g_0(t) & \text{ if } t\in \sT^+ \\
g^{-1}_0(1/t)& \text{ if } t\in \sT\setminus \sT^+
 \end{array}
 \right. ,
 \end{equation*}
is a matching function.

   \item If $g_1, g_2$ are matching functions, then the product $g=g_1 g_2$ is
also a matching function.

    \item If $g$ is a matching function, then for any $a\in L^\infty$
    the duo $(a,ag)$ is a matching pair.

    \item Any matching pair $(a,b)$, $a\in GL^\infty$ can be represented in the form
    $(a,ag)$, where $g=\widetilde{a}\widetilde{b}^{-1}$ is a
    matching function.
 \end{enumerate}

In this section we discuss the Basor-Ehrhardt approach to the study of  defect numbers of the operators $T(a)+H(b)\in \cF(H^p(\sT))$ if $a$ and $b$ are piecewise continuous functions satisfying the condition \eqref{eqn1c}. Then we present an explicit criterion for the Fredholmness of such operators.
Recall that the circle $\sT$ is counterclockwise oriented and $f\in PC$ if and only if for any $t\in\sT$, the one-sided limits
\begin{equation*}
f^{\pm}(t):=\lim_{\ve\to\pm0} f(te^{i\ve})
\end{equation*}
exist. Without loss of generality we assume that $a,b\in GPC$.

\begin{thm}[Basor \& Ehrhardt \cite{BE:2013}]\label{thmBE}
Let $a,b\in GPC$ form a matching pair and let $(c,d)$ be the subordinated pair. The operator $T(a)+H(b)$ is Fredholm on the space $H^p$ if and only if the following conditions hold:
 \begin{align}
 & \frac{1}{2\pi}\arg c^-(1)\notin\left \{\frac{1}{2}+\frac{1}{2p}+\sZ \right \}, \; \frac{1}{2\pi}\arg \widetilde{d}^-(1)\notin\left \{\frac{1}{2}+\frac{1}{2q}+\sZ\right \}, \label{eqn1e}\\
   & \frac{1}{2\pi}\arg c^-(-1)\notin\left \{\frac{1}{2p}+\sZ\right\}, \; \frac{1}{2\pi}\arg \widetilde{d}^-(-1)\notin\left \{\frac{1}{2q}+\sZ\right \}, \label{eqn1f}\\
    & \frac{1}{2\pi}\arg \left ( \frac{c^-(\tau)}{c^+(\tau)} \right )\!\notin\!\left \{\frac{1}{p}\!+\!\sZ\right \}, \; \frac{1}{2\pi}\arg \left (\frac{(\widetilde{d})^-(\tau)}{(\widetilde{d})^+(\tau)}\right)\!\notin\!\left \{\frac{1}{q}\!+\!\sZ\right \},
  \; \forall \tau\in \sT^+.\!\! \label{eqn1g}
 \end{align}
 \end{thm}
The definition of $d$ in \cite{BE:2013} differs from the one used here. In fact, $d$ in \cite{BE:2013} corresponds to $\widetilde{d}$ used here. We keep this notation for sake of easy comparability of the results.

Theorem \ref{thmBE} already shows the exceptional role of the endpoints $+1$ and $-1$ of the upper semicircle $\sT^+$. To establish the index formula  mentioned in \cite{BE:2013}, a more geometric interpretation of the Fredholm conditions \eqref{eqn1e}-\eqref{eqn1g} is needed. Here there are a few details from \cite{BE:2013}.

For $z_1,z_2\in \sC$ and $\theta\in (0,1)$, we consider the open arc $\cA(z_1,z_2,\theta)$ connecting the points $z_1$ and $z_2$ and defined by
 \begin{equation*}
\cA(z_1,z_2,\theta):=\left \{ z\in \sC\setminus \{z_1,z_2\}: \frac{1}{2\pi}\arg\left (\frac{z-z_1}{z-z_2}\right ) \in \left \{\theta +\sZ \right \}\right \}.
 \end{equation*}
For $\theta=1/2$ this arc becomes a line segment and if $z_1=z_2$ it is an empty set. Assuming that $a,b\in GPC$, $a\widetilde{a}=b \widetilde{b}$ and using the auxiliary functions $c$ and $\widetilde{d}$, one can show that the conditions  \eqref{eqn1e}-\eqref{eqn1g} mean that any of the arcs
 \begin{align*}
 & \cA\left (1,c^+(1);\frac{1}{2}+\frac{1}{2p}\right ), \; \cA\left (c^-(\tau), c^+(\tau);\frac{1}{p}\right ),\; \cA\left (c^-(-1),1;\frac{1}{2p}\right ), \\
  & \cA\left (1,(\widetilde{d})^+(1);\frac{1}{2}+\frac{1}{2q}\right ), \; \cA\left ((\widetilde{d})^-(\tau), (\widetilde{d})^+(\tau);\frac{1}{q}\right ),\; \cA\left ((\widetilde{d})^-(-1),1;\frac{1}{2q}\right ),
  \end{align*}
where $\tau\in \sT^+$, does not cross the origin. It is clear that one has to take into account only the jump discontinuity points since $c, \widetilde{d}\in GPC$ and consequently $c^{\pm}(\tau)\neq0$ and $(\widetilde{d})^{\pm}(\tau)\neq0$.

The functions $c$ and $\widetilde{d}$ satisfy the condition $c \widetilde{c}=d \widetilde{d}=1$, so that they are effectively defined by their values on $\sT^+$ only. If we let $\tau$ run along $\sT^+$ from $\tau=1$ to $\tau=-1$, the image of the function $c$ is the curve with possible jump discontinuities, starting at the point $c^+(1)$ and terminating at $c^-(-1)$.
We now add the arcs $\cA(c^-(\tau),c^+(\tau);1/p)$ to any discontinuity point $\tau$ of $c$ located on $\sT^+$. Besides, if necessary we also add the arcs $\cA(1,c^+(1);1/2+1/(2p))$ and $\cA(c^-(-1),1;1/(2p))$ connecting the endpoints $c^+(1)$ and $c^-(-1)$ with the point $\tau=1$, respectively. That way, we obtain a closed oriented curve.  If the operator $T(a)+H(b)$ is Fredholm, the curve does not cross the origin and we consider its winding number $\wind (c^{\#,p})\in \sZ$. Similar constructions lead to the curve $\widetilde{d}^{\#,q}$ with a winding number $\wind (\widetilde{d}^{\#,q})\in \sZ$. Now we can conclude that $T(a)+H(b)\in \cF(H^p)$ if and only if the origin does not belong to the curve $c^{\#,p}$ or $\widetilde{d}^{\#,q}$.

\begin{thm}[Basor \& Ehrhardt \cite{BE:2013}]\label{thmBE1}
Assume that $a,b\in GPC$ form a matching pair with the subordinated pair $(c,d)$ such that the conditions \eqref{eqn1e}-\eqref{eqn1g} hold. Then $T(a)+H(b)\in \cF(H^p)$ with the Fredholm index
\begin{equation*}
\ind (T(a)+H(b))=\wind (\widetilde{d}^{\#,q})-\wind (c^{\#,p}).
\end{equation*}
 \end{thm}
For what follows we need a definition.
\begin{defn}[Basor \& Ehrhardt \cite{BE:2013}]
A matching pair $(a,b)$, $a,b\in L^\infty$ with the subordinated pair $(c,d)$ satisfies the basic factorization condition in $H^p$ if $c$ and $d$ admit the factorization of the form
 \begin{align}
& c(t)=c_+(t)t^{2n}c_+^{-1}(t^{-1}), \quad n\in\sZ,\label{eqn1m_1}\\
& \widetilde{d}(t)=(\widetilde{d})_+(t)t^{2m}(\widetilde{d})_+^{-1}(t^{-1}),\quad m\in\sZ \label{eqn1n_1}
 \end{align}
 and
 \begin{align*}
 &(1+t)c_+(t)\in H^q,\quad (1-t)c_+^{-1}(t)\in H^p , \\%\label{eqn1m_2}  \\
    & (1+t)(\widetilde{d})_+(t)\in H^p,\quad (1-t)(\widetilde{d})_+^{-1}(t)\in H^q .\\% \label{eqn1n_2}
 \end{align*}
 \end{defn}
Note that the indices $m,n$ are uniquely determined and the functions $c_+$ and $d_+$ are also unique up to a multiplicative constant. The representations \eqref{eqn1m_1}, and \eqref{eqn1n_1} are called  antisymmetric factorization of the functions $c$ and $\widetilde{d}$, respectively.
 \begin{thm}[Basor \& Ehrhardt \cite{BE:2013}]\label{thmBE2}
 If $(a,b)$, $a,b\in PC$ is a matching pair and $T(a)+H(b)\in\cF(H^p)$, then $(a,b)$ satisfies the basic factorization condition.
 \end{thm}
The next result allows to determine the defect numbers of the operators $T(a)+H(b)$ in certain situations.

\begin{thm}[Basor \& Ehrhardt \cite{BE:2013}]\label{thmBE3}
Assume that the matching pair\linebreak $(a,b)$, $a,b\in L^\infty$ satisfies the basic factorization condition in $H^p$ with $m,n\in \sZ$. If $T(a)+H(b)\in \cF(H^p)$, then
\begin{align*}
\dim\ker (T(a)+H(b))&=\left \{
 \begin{aligned}
  & 0  & \quad \text{if}\quad n>0, m\leq0,\\
  & -n &\quad \text{if}\quad n\leq0, m\leq0,\\
    & m-n & \quad \text{if}\quad n\leq 0, m>0,\\
   &\dim\ker A_{n,m} &\quad \text{if}\quad n>0, m>0,
 \end{aligned}
 \right .\\
 \dim\ker (T(a)+H(b))^*&=\left \{
 \begin{aligned}
  & 0  &  \text{if}\quad m>0, n\leq0,\\
  & -m & \text{if}\quad m\leq0, n\leq0,\\
    & n-m &  \text{if}\quad m\leq 0, n>0,\\
   &\dim\ker (A_{n,m})^\intercal & \text{if}\quad m>0, n>0,
 \end{aligned}
 \right .
 \end{align*}
Therein, in case $n>0$, $m>0$,
 \begin{equation*}
A_{n,m}:=\left [ \rho_{i-j} + \rho_{i+j}\right ]_{i=0,j=0}^{n-1,m-1}
 \end{equation*}
 and
 \begin{equation*}
\rho(t):=t^{-m-n}(1+t)(1+t^{-1})\widetilde{c}_+\widetilde{d}_+b^{-1}\in L^1(\sT).
 \end{equation*}
 In particular, the Fredholm index of $T(a)+H(b)$ is $m-n$.
 \end{thm}

Let us now briefly discuss the notion of the adjoint operator used in \cite{BE:2013}. If we identify the dual to $H^p$ with $H^q$ via the mapping $g\in H^q\mapsto <g,\cdot >\in (H^p)^\prime$,
\begin{equation*}
<g,f >:=\frac{1}{2\pi}\int_0^{2\pi}g(e^{-i\theta})f(e^{i\theta})\,d\theta,
\end{equation*}
then the adjoint operator to $T(a)+H(b)$ has the form
$T(\widetilde{a})+H(b)$, so that
\begin{equation*}
 \dim\ker (T(a)+H(b))^*=\dim\ker (T(\widetilde{a})+H(b)).
\end{equation*}
A natural question to ask is whether the Fredholmness of $T(a)+H(b)$, $a,b\in L^\infty$ implies the existence of the antisymmetric factorizations of the functions $c$ and $\widetilde{d}$. This problem has been also discussed in \cite{BE:2013} and an example there shows that this is not true in general.

The theory above has been used to examine the operators \eqref{eqn1d} and also the operators $I+H(b)$ with matching functions $b$. We are not going to discuss this approach here. However, in what follows, we present a simple method to handle these operators.

\subsection{Classical approach revisited.\label{ssec2.3}}

Now we turn attention to another method based on the classical approach and recently developed by  the authors of this paper.
Let us start with  a special factorization of the operator $T(V(a,b))$ for generating functions $a$ and $b$ constituting a matching pair. If
$a\in GL^\infty(\sT)$, $b\in L^\infty(\sT)$ satisfy the matching condition \eqref{eqn1c}, the matrix function $V(a,b)$ in \eqref{eqTV} has the form
    \begin{equation}\label{Eq5}
V(a,b)=\left (
 \begin{array}{cc}
 0 & d\\
-c & \widetilde{a}^{-1}
 \end{array}
\right ),\; c=\frac{a}{b}, \; d=\frac{a}{\widetilde{b}}.
    \end{equation}
It follows that the corresponding block Toeplitz operator $T(V(a,b))$ with the generating matrix-function \eqref{Eq5} can be represented in the form   \begin{align}\label{Eq6}
& \quad T(V(a,b))=
\left (
 \begin{array}{cc}
 0 & T(d)\\
-T(c) & T(\widetilde{a}^{-1})
 \end{array}
\right ) \nn\\
&  =
\left (
 \begin{array}{cc}
 - T(d)& 0\\
0 & I
 \end{array}
\right )
\left (
 \begin{array}{cc}
 0& -I\\
I & T(\widetilde{a}^{-1})
 \end{array}
\right )
\left (
 \begin{array}{cc}
 - T(c)& 0\\
0 & I
 \end{array}
\right ),
 \end{align}
with the invertible  operator
\begin{equation*}
\left (
 \begin{array}{cc}
 0& -I\\
I & T(\widetilde{a}^{-1})
 \end{array}
\right ) : H^p\times H^p \to H^p\times H^p .
\end{equation*}
This representations turns out to be extremely useful in the study of Toeplitz plus Hankel operators. To show this we start with the Coburn-Simonenko theorem. For Toeplitz operators this theorem indicates that for any non-zero $a\in L^\infty(\sT)$ one has
\begin{equation*}
\min \{ \dim\ker T(a),\dim\coker T(a) \}=0.
\end{equation*}
It follows that Fredholm Toeplitz operators with  index zero are invertible. However, in general for block-Toeplitz and for Toeplitz plus Hankel operators, Coburn-Simonenko theorem is not true. This causes serious difficulties when studying the invertibility of the operators involved. Nevertheless, the following theorem holds.
 \begin{thm}\label{thm2.9}
Let $a\in GL^\infty(\sT)$ and $A$ refer to one of the operators
$T(a)+H(a t),\; T(a)-H(a t^{-1}),\;T(a)\pm H(a)$. Then
\begin{equation*}
\min\{\dim\ker A, \dim\coker A\}=0.
\end{equation*}
  \end{thm}
  \textit{Proof}.
For $a\in PC(\sT)$ this result goes back to Basor and Ehrhardt  \cite{BE:1999, BE:2004, E:2004h}  with involved proofs. However, there is an extremely simple proof --- cf.~\cite{DS:2014}, based on the representation \eqref{eqTV} and valid for generating functions $a\in GL^\infty$. We would like to sketch this proof here. Thus one of the consequences of the Eq.~\eqref{eqTV} is that there is an isomorphism between the kernels of the operators $T(V(a,b))$ and $\diag (T(a)+H(b),T(a)-H(b))$. Let us start with the operators $T(a)\pm H(a)$. The corresponding subordinated pairs $(c,d)$ have the form $(\pm1, a \widetilde{a}^{-1})$, the third operator in the Eq.~\eqref{Eq6} is $\diag(\mp I,I)$, so that it does not influence the kernel and the image of $T(V(a,\pm a))$. Considering the two remaining operators in \eqref{Eq6}, we note that the Coburn-Simonenko theorem is valid for the block Toeplitz operator $T(V(a,\pm a))$ and hence for $T(a)\pm H(a)$.

Consider now the operators $T(a)+H(at)$. The duo  $(a,at)$ is a matching pair with the subordinated pair $(t^{-1},d)$, $d=a
\widetilde{a}^{-1}t$. The operator $T(t^{-1})$ is Fredholm and $\ind T(t^{-1})=1$, so that $\im T(t^{-1})=H^p$. Besides, a direct check shows that the function $e:=e(t)=1$, $t\in \sT$ belongs to the kernels of both operator $T(t^{-1})$ and  $T(a)-H(at)$. Assuming that $\dim\ker T(d)>0$ and using  Coburn-Simonenko theorem for Toeplitz operators, we note that  $\im T(d)$  is dense in $H^p$. On the other hand, the factorization \eqref{Eq6} and the~Eq.~\eqref{eqTV} yield that both spaces $\im T(V(a,a t))$ and $\im\diag (T(a)+H(at),T(a)-H(a t))$ are dense in $H^p\times H^p$.
Hence,
\begin{equation*}
\coker (T(a)+H(at))=\coker (T(a)-H(at)=\{0\}.
\end{equation*}
Passing to the case $\dim\ker T(d)=0$, we first note that
\begin{equation*}
 1=\dim\ker T(t^{-1}))=\dim\ker\diag((T(a)+H(at),T(a)-H(a t)),
\end{equation*}
and since the kernel of the operator $T(a)-H(a t)$ contains the function $e(t)=1$, it follows that
\begin{equation*}
\ker (T(a)+H(a t))=\{0\}.
\end{equation*}
Thus if $\dim \ker T(a\widetilde{a}^{-1}t)>0$, then $\coker (T(a)+H(at))=\{0\}$, otherwise $\ker (T(a)+H(a t))=\{0\}$ and the Coburn-Simonenko theorem is proved for the operators $T(a)+H(a t)$.
The operators  $T(a)-H(a t^{-1})$ can be considered analogously \cite{DS:2014}.
  \qed

 \vspace{3mm}

Theorem~\ref{thm2.9} can be extended in a few directions --- cf.~Proposition~\ref{prop3.9} and Corollary~\ref{cor3.10} below. It is also valid for Toeplitz plus Hankel operators acting on $l^p$-spaces \cite{DS:2014b} and for Wiener-Hopf plus Hankel integral operators acting on $L^p(\sR^+)$-spaces \cite{DS:2018}.

\section{Kernel Representations.\label{sec3}}

As was already mentioned, the Eq.~\eqref{eqTV}  is of limited use in studying the Fredholmness of the operators  $T(a)+ H(b)$. However, if the generating functions $a$ and $b$ satisfy the matching condition, the representation \eqref{eqTV} allows to determine defect numbers and obtain efficient representations for the kernels and cokernels of Toeplitz plus Hankel operators. In order to present the method,  we recall relevant results for operators acting on Hardy spaces --- cf.~\cite{DS:2014}. Let us start with the connections between the kernels of Toeplitz plus/minus operators and the kernels of the corresponding block Toeplitz operators. The following lemma is a direct consequence of the Eq.~\eqref{eqTV} and is valid even if $a$ and $b$ do not constitute a matching pair.
\begin{lem}\label{l1}
If $a\in GL^\infty$, $b\in L^\infty$ and the operators
$T(a)\pm H(b)$ are considered on the space $H^p$, $1<p<\infty$, then the following relations hold:
  \begin{itemize}
    \item If $(\vp,\psi)^T\in \ker T(V(a,b))$, then
  \begin{align}\label{eqKer1}
     (\Phi, \Psi)^T= &
     \frac{1}{2} (\vp-JQc\vp+JQ\widetilde{a}^{-1}\psi,
      \vp+JQc\vp-JQ\widetilde{a}^{-1}\psi)^T
      \\
     &\in
     \ker \diag(T(a)+H(b),T(a)-H(b)). &\nn
\end{align}

    \item If $(\Phi, \Psi)^T\in \ker\diag (T(a)+H(b),T(a)-H(b)
 )$, then
\begin{equation}\label{eqKer2}
(\Phi+\Psi, P(\widetilde{b}(\Phi + \Psi)
+\widetilde{a}JP(\Phi-\Psi))^T\in \ker T(V(a,b)).
\end{equation}
 \end{itemize}
Moreover, the operators
 \begin{align*}
   E_1:  \ker T(V(a,b)) \to \ker\diag(T(a) + H(b),T(a) - H(b)),\\
       E_2 :  \ker\diag(T(a) + H(b),T(a) - H(b)) \to \ker T(V(a,b)),
\end{align*}
defined, respectively, by the Eqs.~\eqref{eqKer1},
\eqref{eqKer2} are mutually inverses to each other.
 \end{lem}

 \subsection{Subordinated operators and kernels of $T(a)+H(b)$.}

Thus, if the kernel of the operator $T(V(a,b))$ is known,  the kernels of both operators $T(a)+H(b)$ and $T(a)-H(b)$ can be also determined. However, the kernels of the operators $T(V(a,b))$ are known only for a few special classes of matrices $V(a,b)$, and in the case of general generating functions $a,b\in L^\infty$ the kernel $T(V(a,b))$ is not available. The problem becomes more manageable if $a$ and $b$ form a matching pair. In this case, $V(a,b)$ is a triangular matrix --- cf.~\eqref{Eq5} and the subordinated functions $c$ and $d$  satisfy the equations
\begin{equation*}
c \widetilde{c}=1, \quad d\widetilde{d}=1.
\end{equation*}
Moreover, it follows from the Eq.~\eqref{Eq6} that for any function $\vp\in \ker T(c)$, the vector $(\vp,0)^\intercal$ belongs to the kernel of the operator $T(V(a,b))$ and the first assertion in Lemma~\ref{l1} shows that
 \begin{equation}\label{eq5}
  \begin{aligned}
 \vp-JQcP\vp \in \ker(T(a)+H(b)),\\
 \vp+JQcP\vp\in \ker(T(a)-H(b)).
\end{aligned}
 \end{equation}
Another remarkable fact is  that  both functions $\vp+JQcP\vp$ and $\vp-JQcP\vp$ also belong to the kernel of the operator $T(c)$. In order to  show this, we need an auxiliary result.
 \begin{lem}\label{l3}
 Let $g\in L^\infty$ satisfy the relation $g \widetilde{g}=1$ and $f\in \ker T(g)$. Then
 \begin{equation}\label{ker1}
JQgPf\in \ker T(g)\quad\text{and}\quad  (JQgP)^2f=f.
 \end{equation}
 \end{lem}
\textit{Proof}.
We only check the first relation \eqref{ker1}).  Thus
  \begin{equation*}
  T(g)(JQgPf)=PgPJQgPf
  =JQ\widetilde{g}QgPf=JQ\widetilde{g}gPf-JQ\widetilde{g}PgPf=0,
\end{equation*}
so that $JQgPf\in \ker T(g)$.
 \qed
 
 \vspace{3mm}
 
Considering the operators  $\mathbf{P}_g^{\pm}:\ker T(g)\to \ker T(g)$,
\begin{equation*}
\mathbf{P}_g^{\pm}:=\frac{1}{2}(I\pm JQgP)\Big|_{\ker T(g)},
\end{equation*}
we observe that according to Lemma \ref{l3}, one has $(\mathbf{P}_g^{\pm})^2=\mathbf{P}_g^{\pm}$. Therefore, $\mathbf{P}_g^{\pm}$ are complementary projections. This property and the Eqs.~\eqref{eq5} lead to the following conclusion.
 \begin{cor}\label{cor3.3}
If $(c,d)$ is the subordinated pair for the matching pair $(a,b)\in
L^\infty\times L^\infty$, then
  \begin{align*}
   &\ker T(c)=\im \mathbf{P}_c^{-}\dotplus\im\mathbf{P}_c^{+},\\
  & \im \mathbf{P}_c^{-}\subset \ker(T(a) +H(b)) ,\\
   & \im \mathbf{P}_c^{+}\subset \ker(T(a)-H(b)).
\end{align*}
 \end{cor}
Corollary~\ref{cor3.3} shows the impact of the operator $T(c)$ on the kernels
of $T(a)+H(b)$ and $T(a)-H(b)$. The  impact of the operator $T(d)$ on $\ker(T(a)\pm H(b))$ is much more involved.  Thus if $T(c)$ is right-invertible and $s\in \ker T(d)$, then
\begin{equation*}
(T_r^{-1}(c)
T(\widetilde{a}^{-1})s,s)^\intercal \in \ker T(V(a,b)),
\end{equation*}
where $T_r^{-1}(c)$ is one of the right-inverses for $T(c)$.

We now can obtain the representation of the kernel of $T(V(a,b))$.
\begin{lem}
Let $(a,b)$ be a matching pair such that
the operator $T(c)$ is in\-ver\-tible from the right. Then
\begin{equation*}
\ker T(V(a,b))=\Omega(c) \dotplus \widehat{\Omega}(d)
\end{equation*}
 where
 \begin{align*}
\Omega(c) &:=\left \{   (\vp,0)^T:\vp\in \ker T(c)\right\},\\[-0.5ex]
\widehat{\Omega}(d) &:=\left \{ (T_r^{-1}(c)
T(\widetilde{a}^{-1})s,s)^T:s\in \ker T(d) \right \}.
\end{align*}
 \end{lem}

Taking into account this representation and using the first assertion in  Lemma~\ref{l1}, we obtain that if $s\in \ker T(d)$, then
\begin{equation}\label{eqnphi}
 \begin{aligned}
   2 \vp_{\pm}(s)=T_r^{-1}(c)T(\widetilde{a}^{-1}) s&\mp JQcP
   T_r^{-1}(c)T(\widetilde{a}^{-1})s \pm JQ \widetilde{a}^{-1} s\\
  & \quad \in \ker(T(a)\pm H(b))
    \end{aligned}
\end{equation}
The operators $\vp_{\pm}$  can be referred as transition operators, since they transform the kernel of $T(d)$ into the kernels of the operators $T(a)\pm H(b)$. An important property of these operators $\vp_\pm$ is expressed by the following lemma.

\begin{lem}\label{l2}
Let $(c,d)$ be the subordinated pair for a matching pair $(a,b)\in
L^\infty \times L^\infty$. If the operator $T(c)$ is
right-invertible, then for every $s\in \ker T(d)$ the following
relations
\begin{align*}
   (P \widetilde{b} P +P\widetilde{a}JP)\vp_+(s)& = \mathbf{P}_d^+(s), \\
(P \widetilde{b} P -P\widetilde{a}JP)\vp_-(s)& = \mathbf{P}_d^-(s),
\end{align*}
hold. Thus the transition operators $\vp_+$ and $\vp_-$ are injections on the
spaces $\im \mathbf{P}_d^+$ and $\im \mathbf{P}_d^-$, respectively.
 \end{lem}

Using Lemmas \ref{l1}-\ref{l2}, one can obtain a complete description for the kernel of the operator $T(a)+H(b)$ if $(a,b)$ is a matching pair and $T(c)$ is right-invertible.

\begin{prop}[VD \& BS \cite{DS:2014}]\label{prop3.6}
 Let $(c,d)$ be the subordinated pair for the matching pair $(a,b)\in L^\infty \times L^\infty$.
 If the operator $T(c)$ is right-invertible, then the kernels of the operators $T(a)\pm H(b)$ can be represented in the form
 \begin{align*}
   \ker(T(a)+H(b))& =\vp_+(\im \mathbf{P}_d^+) \dotplus\im \mathbf{P}_c^-, \\
\ker(T(a)-H(b))& =\vp_-(\im \mathbf{P}_d^-) \dotplus\im
\mathbf{P}_c^+.
\end{align*}
 \end{prop}
 \begin{rem}\label{rem3.7}
 It was shown in \cite{DS:2014} that  the operators $\vp_+$ and $\vp_-$ send the elements of the spaces $\im \mathbf{P}_d^-$ and $\im \mathbf{P}_d^+$ into $\im \mathbf{P}_c^-$ and $\im \mathbf{P}_c^+$, respectively. Therefore,
 \begin{equation*}
\vp_+: \im \mathbf{P}_d^-\to \im \mathbf{P}_c^-, \quad \vp_-: \im \mathbf{P}_d^+ \to \im \mathbf{P}_c^+
 \end{equation*}
 are well-defined linear operators. If $\im \mathbf{P}_c^-=\{0\}$ ($\im \mathbf{P}_c^+=\{0\}$), then $\vp_+(s_-)=0$ for all $s_-\in \im \mathbf{P}_d^-$ ($\vp_-(s_+)=0$ for all $s_+\in \im \mathbf{P}_d^+$), which yields
 \begin{align*}
  & \vp_-(s_-)= T_r^{-1}(c) T(\widetilde{a}^{-1})s_- , \quad s_-\in \im \mathbf{P}_d^-,\\
    & \left (\vp_+(s_+)= T_r^{-1}(c) T(\widetilde{a}^{-1})s_+ , \quad s_+\in \im \mathbf{P}_d^+\right ).
  \end{align*}
 \emph{Proof}. Indeed, assume for instance that $s_-\in \im \mathbf{P}_d^-$. Then $\vp_+(s_-)=0$ leads to
 \begin{equation*}
   0=2\vp_+(s_-)=T_r^{-1}(c) T(\widetilde{a}^{-1})s_- - JQcPT_r^{-1}(c)T(\widetilde{a}^{-1})s_- +JQ\widetilde{a}^{-1}s_-.
   \end{equation*}
 Hence,
  \begin{equation*}
   T_r^{-1}(c) T(\widetilde{a}^{-1})s_- = JQcPT_r^{-1}(c)T(\widetilde{a}^{-1})s_- -JQ\widetilde{a}^{-1}s_-.
   \end{equation*}
 Thus
 \begin{align*}
 2 \vp_-(s_-)&=T_r^{-1}(c)T(\widetilde{a}^{-1}) s_- + JQcP
   T_r^{-1}(c)T(\widetilde{a}^{-1})s_- - JQ \widetilde{a}^{-1} s_- \\ & =2T_r^{-1}(c)T(\widetilde{a}^{-1}) s_-
 \end{align*}
 and the claim follows. \qed
 
 \vspace{3mm}

 These representations of the transition operators are simpler than \eqref{eqnphi} and it would be interesting to find out which conditions ensure that the restrictions $\vp_+\left |_{\im \mathbf{P}_d^-}\right.$ and $\vp_-\left |_{\im \mathbf{P}_d^+}\right.$ become zero functions.   \end{rem}

Thus in order to obtain an efficient description of the spaces $\ker(T(a)+H(b))$
and  $\ker(T(a)-H(b))$, one has to characterize the projections $\mathbf{P}_c^{\pm}$ and $\mathbf{P}_d^{\pm}$ first.  Such a characterization
can be derived from the Wiener-Hopf factorization \eqref{fac} of the subordinated functions $c$ and $d$. The Wiener-Hopf factorization of these functions can be described in more details, which yields a very comprehensive representation of the kernels of $T(c)$, $T(d)$  and the related projections $\mathbf{P}_c^{\pm},\mathbf{P}_d^{\pm}$.

We first consider related constructions for a matching function $g$ such that  the operator $T(g)$ is Fredholm on $H^p$ with the index $n$. One can show that under the condition $g_-(\infty)=1$, the factorization \eqref{fac} takes the form
 \begin{equation}\label{WHF}
g=\boldsymbol\sigma (g) g_+ t^{-n}\widetilde{g}_+^{-1},
 \end{equation}
where $\boldsymbol\sigma (g)=g_+(0)=\pm1$ and $g_-=\boldsymbol\sigma (g)\widetilde{g}_+^{-1}$. The term $\boldsymbol\sigma (g)$ is called factorization signature. The finding of $\boldsymbol\sigma (g)$ is  a non-trivial problem but if $T(g)$ is invertible and $g$ is continuous at $t=1$ or $t=-1$, then   $\boldsymbol\sigma (g)=g(1)$ or $\boldsymbol\sigma (g)=g(-1)$, respectively. For piecewise continuous functions $g$, the term $\boldsymbol\sigma (g)$ can be also determined.

Notice that  $T(a)-H(b)$ can be also written as Toeplitz plus Hankel operator $T(a)+H(-b)$.  Thereby, the duo $(a,-b)$ is again a matching pair with the subordinated pair $(-c,-d)$ and for the factorization signatures we have $\boldsymbol\sigma (-c)=-\boldsymbol\sigma (c)$, $\boldsymbol\sigma (-d)=-\boldsymbol\sigma (d)$. This observation shows that we can restrict ourselves to the study of Toeplitz plus Hankel operators. Nevertheless, in some cases it is preferable to consider the operator $T(a)-H(b)$ too. But then the leading role still belongs to the operator $T(a)+H(b)$ since the notions of matching pair, subordinated pair and factorization signature are associated with this operator.

Let $g$ stand for the subordinated function $c$ or $d$, so that $g \widetilde{g}=1$. If $T(g)$ is Fredholm, then the factorization  \eqref{WHF} exists with a function $g_+$ satisfying the conditions for factorization \eqref{fac}  and we can describe the spaces $\im \mathbf{P}_g^{\pm}$. This description depends on the evenness of the index of $T(g)$.

\begin{thm}[VD \& BS \cite{DS:2014}]
Assume that $g$ is a matching function, the operator $T(g)$ is Fredholm, $\ind T(g)=n\geq0$ and $g_+$ is the plus factor in the Wiener-Hopf factorization \eqref{WHF} of $g$ in $H^p$. Then
\begin{itemize}
\item For $n=2r$, $r\in \sN$,  the systems of
functions
\begin{equation*}
\cB_{\pm}(g):= \{g_+^{-1}(t^{r-k-1}\pm \boldsymbol\sigma
(g)t^{r+k}):k=0,1,\cdots, r-1 \},
    \end{equation*}
    form bases in the spaces $\im \mathbf{P}_g^{\pm}$ and $\dim\ker \mathbf{P}_g^{\pm}=r$.

\item For $n=2r+1$, $r\in \sZ_+$,  the systems of functions
  \begin{equation*}
\cB_{\pm}(g):= \{g_+^{-1}(t^{r+k}\pm \boldsymbol\sigma
(g)t^{r-k}):k=0,1,\cdots, r \}\setminus \{0\},
  \end{equation*}
form bases in the spaces $\im \mathbf{P}_g^{\pm}$ and $\dim\ker \mathbf{P}_g^{\pm}=r +(1\pm \boldsymbol\sigma (g))/2$.
\end{itemize}
\end{thm}
Thus if $T(c),T(d)$ are Fredholm and $T(c)$ is right-invertible,  Proposition~\ref{prop3.6} provides complete description of the spaces $\ker (T(a)\pm H(b))$. On the other hand, if $T(c)$ is Fredholm but not right-invertible,  the representation
\begin{equation*}
T(a)+H(b)=(T(at^{-n})+H(bt^n)) T(t^n)
\end{equation*}
can be used to study $\ker (T(a)+H(b))$. This is because for any matching pair $(a,b)$ the duo $(at^{-n},b t^n)$ is also a matching pair with the subordinated pair $(ct^{-2n},d)$. A suitable choice of $n$ leads to the right-invertibility of the operator $T(c t^{-2n})$ and we consequently obtain
\begin{equation}\label{BBB}
\ker(T(a)+H(b)) =\ker (T(at^{-n})+H(bt^n)) \cap \im T(t^n).
\end{equation}
The representation \eqref{BBB} has been used in \cite{DS:2014}, to derive the description of the kernels of the operators $T(a)\pm H(b)$. It can be also employed to study one-sided or generalized invertibility of Toeplitz plus Hankel operators and to construct the corresponding one-sided and generalized inverses \cite{DS:2016a,DS:2017,DS:2019a}. In the forthcoming sections, invertibility problems will be  discussed in more details. In this connection,   we note that the (familiar) adjoint operator $(T(a)+H(b))^*$ can be identified with the operator $T(\overline{a})+H(\overline{\widetilde{b}})$ acting on the space $H^q$, $1/p+1/q=1$. It is easily seen that for any matching pair $(a,b)$, the duo $(\overline{a},\overline{\widetilde{b}})$ is also a matching pair with the subordinated pair $(\overline{d},\overline{c})$ and $\boldsymbol{\sigma}(\overline{c})=\boldsymbol{\sigma}(c)$, $\boldsymbol{\sigma}(\overline{d})=\boldsymbol{\sigma}(d)$. Therefore,  cokernel description can be determined directly from the previous results for the kernels of Toeplitz plus Hankel operators.

In some cases the approach above allows omitting the condition of
Fredholmness of the operator $T(d)$.  We note a few results, which can be viewed as an extension of Coburn-Simonenko Theorem~\ref{thm2.9}.

  \begin{prop}[VD \& BS \cite{DS:2014}]\label{prop3.9}
Let $(a,b)\in L^\infty \times L^\infty$ be a matching pair with the
subordinated pair $(c,d)$, and let $T(c)$ be a Fredholm operator.
Then:
\begin{enumerate}[label=(\alph*)]
    \item If $\ind T(c)=1$ and $\boldsymbol\sigma(c)=1$, then
    $$\min\{\dim\ker(T(a)+H(b)), \dim\coker (T(a)+H(b))\}=0.$$

    \item If $\ind T(c)=-1$ and $\boldsymbol\sigma(c)=1$, then
        $$\min\{\dim\ker(T(a)-H(b)), \dim\coker (T(a)-H(b))\}=0.$$

    \item If $\ind T(c)=0$, then
     $$\min\{\dim\ker(T(a)\pm H(b)), \dim\coker (T(a)\pm H(b))\}=0.$$   \end{enumerate}
 \end{prop}
 An immediate consequence of Proposition~\ref{prop3.9} concerns the Toeplitz plus Hankel operators of the form $I\pm H(b)$.
 \begin{cor}\label{cor3.10}
Let $b\in L^\infty$ be a matching function such that
$T(\widetilde{b})$ is a Fredholm operator. Then:
\begin{enumerate}[label=(\alph*)]
    \item If \/ $\ind T(\widetilde{b})=1$
    and $\boldsymbol\sigma(\widetilde{b})=1$, then
     $$\min\{\dim\ker(I+ H(b)), \dim\coker (I+ H(b))\}=0.$$

    \item If\/ $\ind T(\widetilde{b})=-1$
    and $\boldsymbol\sigma(\widetilde{b})=1$, then
   $$\min\{\dim\ker(I- H(b)), \dim\coker (I- H(b))\}=0.$$

        \item If \/ $\ind T(\widetilde{b})=0$,
   then
   $$\min\{\dim\ker(I\pm H(b)), \dim\coker (I\pm H(b))\}=0.$$
\end{enumerate}
 \end{cor}

\subsection{Kernels of $T(a)+H(b)$ for piecewise continuous generating  functions.}

If  more information about generating functions is available, then the kernel of the Fredholm operator $T(a)+H(b)\in \cL(H^p)$  can be studied under weaker conditions. Thus for  piecewise continuous functions $a$ and $b$, the assumption that the subordinated operators $T(c),T(d)\in \cL(H^p)$ are Fredholm can be dropped. In order to show this we need a few facts from \cite{DS:2014}.

Let $A$ be an operator defined on all spaces $L^p$ for $1<p<\infty$.
Consider the set
 $$
 A_{F}:=\{p\in(1,\infty) \; \text{such that the
operator} \; A:H^p\to H^p \; \text{is Fredholm}\}.$$

\begin{prop}[\v{S}ne{\u\i}berg {\cite{Sh:1974}}]
The set $A_F$ is open. Moreover, for each connected component
$\gamma\in A_F$, the index of the operator $A:L^p\to L^p$, $p\in
\gamma$ is constant.
 \end{prop}

 This result also holds for operators $A$ acting on the spaces $H^p$, $1<p<\infty$, since any operator $A\colon H^p\to H^p$ can be identified with the operator $AP+Q$ acting already on $L^p$.  For Toeplitz operators the structure of the set $A_F$ can be characterized as follows.
 \begin{prop}[Spitkovski{\u\i} {\cite{Sp:1976}}]\label{p5}
 Let $G$ be an invertible matrix-func\-tion with entries from $PC$
 and let $A:=T(G)$. Then there is an at most countable subset $S_A\subset
 (1,\infty)$ with the only possible accumulation points $t=1$ and
 $t=\infty$ such that $A_F=(1,\infty)\setminus S_A$.
\end{prop}

Clearly, if $G$ is piecewise continuous with only a finite number of discontinuities, then $S_A$ is at most finite. This result can be used to describe the corresponding set $A_F$ for Toeplitz plus Hankel operators with $PC$-generating functions.

\begin{cor}
Let $a,b\in PC$ and  $A:=\diag(T(a)+H(b),T(a)-H(b))\colon$  $H^p\times
H^p \to H^p\times H^p$. Then there is at most countable subset
$S_A\subset (1,\infty)$ with only possible accumulation points
$t=1$ and $t=\infty$ such that $A_F=(1,\infty)\setminus S_A$.
 \end{cor}
 \textit{Proof}.
It follows directly from Proposition \ref{p5} since
 $\diag(T(a)+H(b),T(a)-H(b))$ is Fredholm if and only if so is the
 operator $T(V(a,b))$.
 \qed

 \vspace{3mm}

Thus if $a,b\in PC$ and the operator $T(a)+H(b)$ is Fredholm on
$H^p$,  there is an interval  $(p',p'')$ containing $p$ such
that $T(a)+H(b)$ is Fredholm on all spaces $H^r$, $r\in (p',p'')$
and the index of this operator does not depend on $r$. Moreover,
there is an interval $(p,p_0)\subset (p',p'')$, $p<p_0$ such that
$T(a)-H(b)$ is Fredholm on $H^r$, $r\in (p,p_0)$ and its index does
not depend on $r$. Recalling that for $\nu<s$ the space $H^s$ is continuously embedded into $H^\nu$, we note that the kernel of $T(a)+H(b):H^r\to H^r$ does not depend on $r\in (p',p'')$. The same is also true for $\ker(T(a)-H(b))$, $r\in (p,p_0)$. 

Hence, the kernel of the operator $T(a)+H(b)$ acting on the space $H^p$ coincides with the kernel of this operator acting on the space $H^r$ for $r\in (p,p_0)$ and the latter can be studied by the approach above. Therefore, if $T(a)+H(b)\in \cL(H^p)$ is Fredholm and $a,b\in PC$ form a matching pair,  the kernel of the operator $T(a)+H(b)$ can be described.

\section{Generalized\! Invertibility. One-Sided\! Inver\-tibility. Invertibility.\label{sec4}}

Let $a,b\in L^\infty$ be a matching pair with the subordinated pair $(c,d)$. The pair is called a Fredholm matching pair (with respect to $H^p$) if the operators $T(c),T(d)\in \mathcal{L}(H^p)$ are Fredholm. We  write $T(a)+H(b)\in \cF_{TH}^p(\kappa_1,\kappa_2)$ if $(a,b)$ is a Fredholm matching pair with the subordinated pair $(c,d)$ such that  $\ind T(c)=\kappa_1$, $\ind T(d)=\kappa_2$. It was first noted in \cite{DS:2014a} that if $\kappa_1\geq 0,\kappa_2\geq 0$ or if $\kappa_1\leq 0,\kappa_2\leq 0$, then \eqref{eqTV} and \eqref{Eq6} yield one-sided invertibility of the operator $T(a)+H(b)$. However, if $\kappa_1\kappa_2< 0$, the invertibility issues become more involved. We start this section by considering the generalized invertibility of the operators $T(a)+H(b)$, $a,b\in L^\infty$. Set
\begin{equation*}
B:=\cP V(a,b)\cP+\cQ,
\end{equation*}
where $\cP:=\diag(P,P)$, $\cQ:=\diag(Q,Q)$.

\begin{thm}[VD \& BS \cite{DS:2016a}]\label{thm4.1}
Assume that $(a,b)$ is a matching pair with the subordinated pair
$(c,d)$ and $B$ is generalized invertible operator, which has a
generalized inverse $B^{-1}$ of the form
 \begin{equation}\label{Eq52}
  B^{-1}=
 \left(%
\begin{array}{cc}
 \bA  & \bB \\
   \bD & 0 \\
   \end{array}%
\right) +\cQ,
 \end{equation}
where  $\bA,\bB$ and $\bD$ are operators acting in the space $H^p$.
Then the operator $T(a)+H(b):H^p\to H^p$ is also generalized invertible and one of its generalized inverses has the form
 \begin{align*}
 (T(a)+H(b))^{-1}_g&=-H(\widetilde{c})(\bA(I-H(d))-\bB
H(\widetilde{a}^{-1})) \\
&\quad +H(a^{-1})\bD(I-H(d))+T(a^{-1}).
\end{align*}
\end{thm}
This result can now be used to construct generalized inverses for the operator $T(a)+H(b)$ in the following cases --- cf.~\cite{DS:2016a}:
  \begin{enumerate}[label=(\alph*)]
    \item If $\kappa_1\geq 0$ and $\kappa_2 \geq 0$, then $B$ has
generalized inverse of the form \eqref{Eq52} with
$\bA=T_r^{-1}(c)T(\widetilde{a}^{-1})T_r^{-1}(d)$,
$\bB=-T_r^{-1}(c)$ and $\bD=T_r^{-1}(d)$.

    \item If $\kappa_1\leq0$ and $\kappa_2\leq0$, then $B$ has
generalized inverse of the form \eqref{Eq52} with
$\bA=T_l^{-1}(c)T(\widetilde{a}^{-1})T_l^{-1}(d)$,
$\bB=-T_l^{-1}(c)$ and $\bD=T_l^{-1}(d)$.

    \item If $\kappa_1\geq 0$ and $\kappa_2\leq0$, then $B$ has
generalized inverse of the form \eqref{Eq52} with
$\bA=T_r^{-1}(c)T(\widetilde{a}^{-1})T_l^{-1}(d)$,
$\bB=-T_r^{-1}(c)$ and $\bD=T_l^{-1}(d)$.
 \end{enumerate}
It is clear that in the cases (a) and (b), generalized inverses are, respectively, right and left inverses. We also note that under conditions of assertion (a), a right inverse of $T(a)+H(b)$ can be written in a simpler form --- cf.~\cite{DS:2017}
  \begin{equation}\label{EqRI}
  B:= (I - H(\widetilde{c})) T_r^{-1}(c) T(\widetilde{a}^{-1})T_r^{-1}(d) + H(a^{-1})T_r^{-1}(d).
 \end{equation}
 The proof of this result is straightforward --- i.e. one can use the relations \eqref{Widom} to verify that $(T(a)+H(b))B=I$. On the other hand, under conditions of (b), a simpler representation of the left-inverse of $T(a)+H(b)$ can be derived from \eqref{EqRI} by passing to the adjoint  operator.

 In addition to the cases considered, there is one more situation --- viz.
  \begin{enumerate}
   \item[(d)] $\kappa_1 <0$, $\kappa_2>0$.
   \end{enumerate}
This case is much more involved and we will deal with it later on. At the moment, we focus on invertibility of operators from
$\cF^p_{TH}(\kappa_1,\kappa_2)$, $1<p<\infty$.

 \begin{thm}[VD \& BS \cite{DS:2017}]\label{thm4.2}

 Assume that $T(a)+H(b)\in \cF^p_{TH}(\kappa_1,\kappa_2)$ is invertible. Then:
\begin{enumerate}[label=(\roman*)]

   \item If $\kappa_1\geq \kappa_2$ or $\kappa_1 \kappa_2 \geq 0$, then
       \begin{equation}\label{cond4.4}
|\kappa_1|\leq 1, \quad|\kappa_2|\leq 1.
       \end{equation}

   \item If $\kappa_1<0$ and $\kappa_2>0$, then
   \begin{enumerate}[label=(\alph*)]

     \item If $\kappa_1$ and $\kappa_2$ are even numbers, then $\kappa_2=-\kappa_1$.

    \item If $\kappa_1$ is an odd number and $\kappa_2$ is an even one, then $\kappa_2=-\kappa_1+\boldsymbol\sigma(c)$.

     \item If $\kappa_1$ is an even number and $\kappa_2$ is an odd one, then $\kappa_2=-\kappa_1-\boldsymbol\sigma(d)$.

     \item  If $\kappa_1$ and $\kappa_2$ are odd numbers, then $\kappa_2=-\kappa_1+\boldsymbol\sigma(c)-\boldsymbol\sigma(d)$.
       \end{enumerate}
       \end{enumerate}
       \end{thm}
Our next goal is to obtain sufficient invertibility conditions for the invertibility of the operators from $\cF^p_{TH}(\kappa_1,\kappa_2)$ and to provide effective representations for their inverses. We assume first that $\kappa_1$ and $\kappa_2$ satisfy conditions \eqref{cond4.4}.
 \begin{thm}[VD \& BS \cite{DS:2017}]\label{thm4.3}
Assume that $T(a)+H(b)\in\cF^p_{TH}(\kappa_1,\kappa_2)$ and  the indices of the subordinated operators $T(c)$, $T(d)$ and the factorization signatures of $c$ and $d$ satisfy one of the following conditions:
 \begin{enumerate}[label=(\roman*)]
   \item $\kappa_1=0$, $\kappa_2=0$;
   \item $\kappa_1=1$, $\kappa_2=0$ and $\boldsymbol\sigma (c)=1$;
   \item $\kappa_1= 0$, $\kappa_2=1$ and $\boldsymbol\sigma (d)=-1$;
   \item $\kappa_1=1$, $\kappa_2=1$ and $\boldsymbol\sigma (c)=1$, $\boldsymbol\sigma (d)=-1$;
   \item $\kappa_1=0$, $\kappa_2=-1$ and $\boldsymbol\sigma (d)=1$;
   \item $\kappa_1=-1$, $\kappa_2=0$ and $\boldsymbol\sigma (c)=-1$;
   \item $\kappa_1=-1$, $\kappa_2=-1$ and $\boldsymbol\sigma (c)=-1$, $\boldsymbol\sigma (d)=1$;
   \item $\kappa_1=1$, $\kappa_2=-1$ and $\boldsymbol\sigma (c)=1$, $\boldsymbol\sigma (d)=1$;
   \item $\kappa_1=-1$, $\kappa_2=1$ and $\boldsymbol\sigma (c)=-1$, $\boldsymbol\sigma (d)=-1$.
 \end{enumerate}
Then the operator $T(a)+H(b)$ is invertible.  Moreover, we have:
 \begin{enumerate}

\item Under conditions (i)-(iv), the inverse operator has the form \eqref{EqRI}, whe\-re right-inverses of $T(c)$ or/and $T(d)$ shall be replaced by the corresponding inverses.

    \item Under conditions (v)-(vii), the inverse operator has the form
        \begin{equation*}%\label{Eq61}
   \begin{aligned}
&(T(a)+H(b))^{-1}=-H(\widetilde{c})(T_l^{-1}(c)T(\widetilde{a}^{-1})T_l^{-1}(d)(I-H(d))\\
 &\quad +T_l^{-1}(c)
H(\widetilde{a}^{-1}))+H(a^{-1})T_l^{-1}(d)(I-H(d))+T(a^{-1}).
 \end{aligned}
 \end{equation*}
 \item Under condition (viii), the inverse operator has the form
 \begin{equation*}%\label{Eq61}
   \begin{aligned}
&(T(a)+H(b))^{-1}=-H(\widetilde{c})(T_r^{-1}(c)T(\widetilde{a}^{-1})T_l^{-1}(d)(I-H(d))\\
 &\quad +T_r^{-1}(c)
H(\widetilde{a}^{-1}))+H(a^{-1})T_l^{-1}(d)(I-H(d))+T(a^{-1}).
 \end{aligned}
 \end{equation*}

 \item Under condition (ix), the inverse operator has the form
  \begin{align*}
&\quad (T(a)+H(b))^{-1}=T(t^{-1})(I - c_+^{-1} t Q t^{-1})\\
&\times [(I - H(t^2\widetilde{c}))
T_r^{-1}(t^{-2}c) T(\widetilde{a}^{-1}t^{-1})T_r^{-1}(d) +
H(a^{-1}t)T_r^{-1}(d) ],\nn
 \end{align*}
where $c_+$ is the plus factor in factorization \eqref{WHF} for the function $c$.
 \end{enumerate}
 \end{thm}
Theorem~\ref{thm4.3} is, in fact, the collection of various results from \cite{DS:2017}. On the other hand, conditions (i)-(ix) are not necessary for the invertibility of $T(a)+H(b)$ in the case \eqref{cond4.4}. Thus if $\kappa_1=-1$, $\kappa_2=1$, then the operator $T(a)+H(b)$ can be invertible even if $\boldsymbol\sigma (c)$ and $\boldsymbol\sigma (d)$  do not satisfy condition (ix). This case, however, requires special consideration and has been not yet studied.

Consider now the situation (ix) in more detail. This is a subcase of assertion (ii) in Theorem~\ref{thm4.2} and a closer inspection shows substantial difference from the other cases in Theorem~\ref{thm4.3}. What makes it very special is the presence of the factorization factor of $c$ in the representation of the inverse operator. It is also worth noting that the construction of $(T(a)+H(b))^{-1}$ is more involved and requires the following result.

\begin{lem}[VD \& BS \cite{DS:2017, DS:2019a}]\label{lem_inverse}
Let $C,D$ be operators acting on a Banach space $X$. If $A=CD$ is an invertible operator, then $C$ and $D$ are, respectively, right- and left-invertible operators. Moreover, the operator $D:X\to \im D$ is invertible and if $D_0^{-1}:\im D\to X$ is the corresponding inverse, then the operator $A^{-1}$ can be represented in the form
\begin{equation*}
A^{-1}=D_0^{-1} P_0 C_r^{-1},
\end{equation*}
where $P_0$ is the projection from $X$ onto $\im D$ parallel to $\ker C$ and $C_r^{-1}$ is any right-inverse of $C$.
 \end{lem}

Theorem \ref{thm4.2}(ii) provides necessary conditions of the invertibility of $T(a)+H(b)\in \cF^p_{TH}(\kappa_1,\kappa_2)$ for negative $\kappa_1$ and positive $\kappa_2$. In the next section, we take a closer look at the condition (iia). The related ideas can be seen as a model to study invertibility in cases (iib)-(iid) of Theorem~\ref{thm4.2}.

Now we would like to discuss a few examples.

\begin{ex}
Let us consider the operator $T(a)+H(b)$ in the case where $a=b$.
In this situation $c(t)=1$ and $d(t)=a(t) \widetilde{a}^{-1}(t)$.
Hence, $H(\widetilde{c})=0$, $T(c)=I$ and if $\ind T(d)=0$, then the
operator $T(a)+H(a)$ is also invertible and
  \begin{equation*}
(T(a)+H(a))^{-1}=(T(\widetilde{a}^{-1}) + H(a^{-1}))T^{-1}(a
\widetilde{a}^{-1}).
  \end{equation*}
\end{ex}

\begin{ex}
Similar approach show that the operator $T(a)+H(\widetilde{a})$ is invertible if $T(c)$,  $c(t)=a(t)
\widetilde{a}^{-1}(t)$ is invertible and
  \begin{equation*}
(T(a)+H(\widetilde{a}))^{-1}=(I - H(\widetilde{a}a^{-1})) T^{-1}(a
\widetilde{a}^{-1}) T(\widetilde{a}^{-1}) + H(a^{-1}).
  \end{equation*}
\end{ex}

\begin{ex}
Consider the operator $I+H(\phi_0 b)$, where $b(t)\widetilde{b}(t)=1$ and $\phi_0(t)=t$, $\phi_0(t)=-t^{-1}$ or $\vp_0(t)=\pm1$ for all $t\in\sT$. Assume that the operator $T(\widetilde{b})$ is Fredholm.
Corollary~\ref{cor3.10} indicates that if one of the conditions
\begin{enumerate}[label=(\alph*)]
    \item  $\ind T(\widetilde{b})=0$ and $\vp_0(t)=\pm1$,
    \item $\ind T(\widetilde{b})=0$, $\boldsymbol\sigma(\widetilde{b})=1$ and $\vp_0(t)=t$,
    \item $\ind T(\widetilde{b})=0$,
   $\boldsymbol\sigma(\widetilde{b})=1$ and $\vp_0(t)=-t^{-1}$,
  \end{enumerate}
   holds, then
   \begin{equation*}
 \min\{ \dim \ker(I+H(\vp_0 b)), \dim\coker(I+H(\vp_0 b))\}=0.
\end{equation*}
 Therefore, if  $I+H(\vp_0 b)$ is Fredholm with index zero, then this operator is invertible. However, for $b\in PC$,  the Fredholmness of the operators $T(\widetilde{b})$ and $I+H(\vp_0 b)$ can be studied by Theorems~\ref{thmBE} and~\ref{thmBE1}. It is also possible to construct the inverse operator using the corresponding results on the factorization of $PC$-functions. However, instead of going into details, we would like to observe that if $T(b)$, $b\in L^\infty$ is invertible, then $I+T(\vp_0 b)$ is also invertible under the conditions of Theorem~\ref{thm4.3}(i), (viii), and (ix), respectively. Moreover, the inverse operator can be explicitly constructed.

 Using a distinct method, Basor and Ehrhardt \cite{BE:2017} also proved the invertibility of this operator on $H^2$ under the condition that $T(\widetilde{b}):H^2\to H^2$ is invertible. For the $H^2$-space, the invertibility of $T(b)$  automatically follows from that of $T(\widetilde{b})$. However, if $p\neq 2$, this is not true and the corresponding examples can be found already among operators with piecewise continuous generating functions. It is interesting enough that in each case $\vp_0(t)=\pm1$, $\vp_0(t)=t$ or $\vp_0(t)=-t^{-1}$, the inverse operator can be represented in the form
  \begin{equation*}
  (I+H(\vp_0 b))^{-1}=(T(\widetilde{b})+H(\vp_0))^{-1}(I+H(\vp_0 \widetilde{b}))(T(b)+H(\vp_0))^{-1}.
  \end{equation*}
  \end{ex}
\section{Invertibility \!of Operators\! from\! $\cF^p_{TH}(\!-2n,2n)$,\\ $n \in \sN$.\label{sec5}}
Theorem~\ref{thm4.2}(ii) provides necessary conditions for invertibility of $T(a)+H(b)\in \cF^p_{TH}(\kappa_1, \kappa_2)$ if $\kappa_1<0$ and $\kappa_2>0$. There are four different situations to consider. Here we focus only on the case $T(a)+H(b)\in \cF^p_{TH}(-2n,2n)$, $n\in \sN$,  but the reader can  handle the remaining cases using the ideas below. Let us start with the simplest case $T(a)+H(b)\in \cF^p_{TH}(-2,2)$ and let $(c,d)$ be the subordinated pair for the matching pair $(a,b)\in L^\infty\times L^\infty$. In passing note that if $T(a)+H(b)\in \cF^p_{TH}(-2,2)$, then the adjoint operator $(T(a)+H(b))^*=T(\overline{a})+H(\widetilde{\overline{b}})$ belongs to the set $\cF^q_{TH}(-2,2)$.

According to \eqref{WHF}, the Wiener-Hopf factorization of the function $d$ is
\begin{equation*}
d(t)=\boldsymbol\sigma (d)d_+(t)\, t^{-2}\, \widetilde{d}_+^{-1}(t).
\end{equation*}
It is easily seen that the operator $T(a)+H(b)$ can be represented in the form
\begin{equation}\label{eqn5.1}
T(a)+H(b)=(T(a_1)+H(b_1))T(t),
\end{equation}
where $a_1=at^{-1}$ and $b_1=tb$. The duo $(a_1,b_1)=(at^{-1},bt)$ is a matching pair with the subordinated pair $(c_1,d_1)=(ct^{-2},d)$. Hence, $T(a_1)+H(b_1)\in \cF^p_{TH}(0,2)$ and we note that $T(c_1)$ is invertible. The invertibility of $T(c_1)$ implies that both  projections $\mathbf{P}_{c_1}^+$ are $\mathbf{P}_{c_1}^-$ are the zero operators. According to Remark~\ref{rem3.7}, the functions $\vp_{\pm}$ admit the representations
\begin{equation*}
     \vp_{\pm}(s_{\pm})= T_r^{-1}(c) T(\widetilde{a}^{-1})s_{\pm} , \quad s_{\pm}\in \im \mathbf{P}_d^{\pm}.
   \end{equation*}
so that
\begin{equation*}%\label{eqn6.5}
\ker (T(a_1)\pm H(b_1))=T_r^{-1}(c_1) T(\widetilde{a}_1^{-1}) (\im \mathbf{P}_d^{\pm}).
\end{equation*}
Further, we also note that
\begin{equation*}
\im \mathbf{P}_d^{\pm}=\{\nu d_+^{-1}(1\pm \boldsymbol\sigma(d)t : \nu\in\sC\}
\end{equation*}
is a one-dimensional subspace of $\ker T(d)$.

By $\omega^{a,b,\pm}$ we denote the functions
  \begin{equation}\label{eqn01}
\omega^{a,b,\pm}(t)=T^{-1}(ct^{-2}) T(\widetilde{a}^{-1}t^{-1})(d_+^{-1}(t)\pm  \boldsymbol\sigma (d)  d_+^{-1}(t)t),
\end{equation}
which respectively  belong to the kernels of the operators $T(a_1)\pm H(b_1)$. It is clear that $\omega^{a,b,\pm}$ also depend on $p$.

Representation \eqref{eqn5.1} shows that $T(a)+H(b)$ has trivial kernel if and only if
\begin{equation*}%\label{eqn9}
\ker (T(a_1)+H(b_1))\cap \im T(t) =\{ 0\}.
 \end{equation*}
It is possible  only if
\begin{equation*}
\widehat{\omega}^{a,b,+}_0 \neq 0,
\end{equation*}
where $\widehat{\omega}^{a,b,+}_0$ is the zero Fourier coefficient of the function $\omega^{a,b,+}(t)$. Similar result is valid for the operator $T(a)-H(b)$. Note that if $T(a)+H(b)$ belongs to the set $\cF^p_{TH}(-2,2)$, then so is the operator $T(a)-H(b)$, since $T(a)-H(b)=T(a)+H(-b)$.

 \begin{thm}[VD \& BS \cite{DS:2019a}]\label{thm5.1}
 Let $T(a)+ H(b)\in \cF^p_{TH}(-2,2)$.
 \begin{enumerate}[label=(\alph*)]
\item The operator $T(a)+H(b)$ $(T(a)-H(b))$ is left-invertible if and only if \/  $\widehat{\omega}^{a,b,+}_0 \neq 0$ $(\widehat{\omega}^{a,b,-}_0 \neq 0)$.

 \item The operator $T(a)+H(b)$ $(T(a)-H(b))$ is right-invertible if and only if \/
 $\widehat{\omega}^{\overline{a},\widetilde{\overline{b}},+}_0 \neq 0$ $(\widehat{\omega}^{\overline{a},\widetilde{\overline{b}},-}_0 \neq 0)$.

 \item The operator $T(a)+H(b)$ $(T(a)-H(b))$ is invertible if and only if \/  $\widehat{\omega}^{a,b,+}_0 \neq 0$ and\/ $\widehat{\omega}^{\overline{a},\widetilde{\overline{b}},+}_0 \neq 0$ $(\widehat{\omega}^{a,b,-}_0 \neq 0$ and\/ $\widehat{\omega}^{\overline{a},\widetilde{\overline{b}},-}_0 \neq 0)$.

 \item If $\widehat{\omega}^{a,b,+}_0 \neq 0$ and $\widehat{\omega}^{a,b,-}_0 \neq 0$, then both operators $T(a)+H(b)$ and $T(a)-H(b)$ are invertible.
  \end{enumerate}
 \end{thm}
Let us sketch the proof of Assertion (d). It follows from the representations \eqref{eqTV} and \eqref{Eq6} that
\begin{equation}\label{AAA}
\ind (T(a)+H(b))+\ind (T(a)-H(b))=0.
\end{equation}
By Assertion (a), both operators $T(a)+H(b)$ and $T(a)-H(b)$ are left-invertible. Therefore, $\ind(T(a)+H(b)) \leq 0$, $\ind(T(a)-H(b)) \leq 0$
and taking into account \eqref{AAA}, we obtain that
\begin{equation*}
\ind (T(a)+H(b))=\ind (T(a)-H(b))=0,
\end{equation*}
which yields the invertibility of both operators under consideration.

If an operator $T(a)+H(b)\in \cF^p_{TH}(-2,2)$ is invertible, we can
 construct its inverse by using Lemma~\ref{lem_inverse}. In particular, we have.

\begin{thm}[VD \& BS \cite{DS:2019a}]
If  $T(a)+H(b)\in \cF^p_{TH}(-2,2)$ is invertible, then the inverse operator can be represented in the form
 \begin{align}\label{inverse}
&(T(a)+H(b))^{-1}= T(t^{-1})\nn\\
&\times\left (I-\frac{1}{\widehat{\omega}^{a,b,+}_0}T^{-1}(ct^{-2}) T(\widetilde{a}^{-1}t^{-1})d_+^{-1}(t)(1+ \boldsymbol\sigma (d)t)tQt^{-1} \right)\nn\\
 &\quad\times \left ((I  -  H(t^2\widetilde{c})) T^{-1}(t^{-2}c)
T(\widetilde{a}^{-1}t^{-1})T_r^{-1}(d)  +  H(a^{-1}t)T_r^{-1}(d)\right )  .
 \end{align}
\end{thm}

\begin{ex}
We now consider the operator $A=T(a)+H(t^{-2}a)$, defined by the function

\begin{equation*}%\label{eqn13}
a(t):=(1-\gamma t^{-1})(1-\gamma t)^{-1},
\end{equation*}
where $\gamma$ is a fixed number in the interval $(0,1)$.
\end{ex}
The duo $(a, at^{-2})$ is a matching pair with the subordinated pair $(c,d)=(t^2,a ^{-1} t^2)$ and $A\in \cF^p_{TH}(-2,2)$. The corresponding Wiener-Hopf factorizations of $a$ and $d$ are the same in all $H^p$. More precisely, we have
\begin{alignat*}{2}
 a(t)&=a_+(t) \widetilde{a}_+^{-1}(t), &\quad  a_+(t)=(1-\gamma t)^{-1} ,&\\
 d(t)&=1 \cdot d_+(t) t^{-2} \widetilde{d}_+^{-1}(t),&
d_+(t)=(1-\gamma t)^{-2}.& %, \quad \widetilde{d}_+^{-1}(t)=(1-\gamma t^{-1})^2.
\end{alignat*}
Hence, $\boldsymbol\sigma(d)=1$ and computing the zero Fourier coefficients of the corresponding functions \eqref{eqn01}, we obtain   $$\widehat{w}^{a,at^{-2},+}_0=\widehat{w}^{a,at^{-2},-}_0=\gamma^2-\gamma+1.$$
It is easily seen that for any $\gamma\in (0,1)$ these coefficients are not equal to zero, so that by Theorem~\ref{thm5.1}(c), the operator $T(a)+H(at^{-2})$ is invertible. The corresponding inverse operator, which is constructed according to the representation \eqref{inverse}, has the form
\begin{align}  \nn
&\quad (T(a)+ H(t^{-2}a))^{-1} \\
&= T(t^{-1}) \left (I-\frac{1}{\gamma^2-\gamma+1}P\left ( (1-\gamma t^{-1})(1-\gamma t)(1+t)\right)Qt^{-1} \right ) \nn  \\
& \quad \times \left (T\left ( \frac{(1-\gamma t^{-1})t^{-1}}{1-\gamma t}\right ) +H \left ( \frac{(1-\gamma t)t}{1-\gamma t^{-1}}\right ) \right )\nn\\& \quad\quad \times T((1-\gamma t)^2)T((1-\gamma t^{-1})^{-2})T(t^2).\nn
\end{align}

 Consider now the invertibility of the operators  $T(a)+H(b)$ from the set  $\cF^p_{TH}(-2n,2n)$ for $n$  greater than $1$. Thus we assume that the subordinated operator $T(c)=T(a b^{-1})$ and $T(d)=T(a\widetilde{b}^{-1})$ are Fredholm and
 \begin{equation*}
\ind T(c)=-2n,\quad \ind T(d)=2n.
 \end{equation*}
  Considering the functions
   \begin{equation*}
\omega_{k}^{a,b}(t) :=T^{-1}(ct^{-2n}) T(\widetilde{a}^{-1}t^{-n})d_+^{-1}(t)(t^{n-k-1}+ \boldsymbol\sigma (d) t^{n+k}), \; k=0,1,\cdots, n-1,
 \end{equation*}
we introduce the matrix
\begin{equation*}
W_n(a,b)=(\omega_{jk}^{a,b})_{k,j=0}^{n-1},
\end{equation*}
with the entries
 \begin{equation*}%\label{eqn30}
\omega_{jk}^{a,b}=\frac{1}{2\pi} \int_{0}^{2\pi} \omega_{k}^{a,b}(e^{i \theta})e^{-i j \theta}\,d\theta, \quad j,k=0,1,\cdots, n-1,
 \end{equation*}
and the terms $d_+$ and $\boldsymbol\sigma (d)$ defined by the Wiener-Hopf factorization
\begin{equation*}
d(t)=\boldsymbol\sigma (d)d_+(t)\, t^{-2n}\, \widetilde{d}_+^{-1}(t),
\end{equation*}
with respect to $L^p$. Notice that the functions $\omega_{k}^{a,b}$ form a basis in $\vp_+(\im\mathbf{P}_d^+)$, where $\vp_+: \im\mathbf{P}_d^+ \to \ker(T(a_1)+H(b_1))$ is defined by
\begin{equation*}
\vp_+=T^{-1}(ct^{-2n}) T(\widetilde{a}^{-1}t^{-n})).
\end{equation*}
The invertibility of the operators from $\cF^p_{TH}(-2n,2n)$ is described by the following theorem.

\begin{thm}[VD \& BS \cite{DS:2019a}]\label{thm5.4}
If $T(a)+H(b)\in \cF_{TH}^p(-2n,2n)$, then:
\begin{enumerate}[label=(\alph*)]
  \item $T(a)+H(b)$ is left-invertible if and only if $W_n(a,b)$ is a non-degenerate matrix.
  \item $T(a)+H(b)$ is right-invertible if and only if $W_n(\overline{a},\widetilde{\overline{b}})$ is a non-degenerate matrix.
  \item $T(a)+H(b)$ is invertible if and only if $W_n(a,b)$ and $W_n(\overline{a},\widetilde{\overline{b}})$ are non-degenerate matrices.
  \item If $W_n(a,b)$ and $W_n(a,-b)$ are non-degenerate matrices, then both operators $T(a)+H(b)$ and $T(a)-H(b)$ are invertible.
\end{enumerate}
\end{thm}
 \begin{ex}
 Consider operator $T(a)+H(at^{-2n})$, $n\in\sN$ and
 \begin{equation}\label{leftn}
a(t)=\frac{1-\gamma t^{-1}}{1-\gamma t}, \; \gamma\in(0,1).
 \end{equation}
This operator is invertible and the inverse operator can be constructed.
 \end{ex}
 It was shown in \cite{DS:2019a} that for any $n\in \sN$, the operator \eqref{leftn} is left-invertible in any space $H^p$, $1<p<\infty$. Moreover, since $H(at^{-2n})$ is compact and $\ind T(a)=0$, the operator at hand is even invertible.

 \begin{rem}
 If $m,n\in \sN$ and $m\neq n$, the set $\cF^p_{TH}(-2m,2n)$ does not contain invertible operators, but it still includes one-sided invertible operators.
  \end{rem}

\section{Toeplitz plus Hankel Operators on\\ $l^p$-spaces.\label{sec6}}

A substantial portion of the results presented in Sections~\ref{sec2}-\ref{sec5} can be extended to Toeplitz plus Hankel operators acting on $l^p$-spaces, $1<p<\infty$. Such an extension is highly non-trivial because many tools perfectly working in classical Hardy spaces $H^p$, are not available for operators on $l^p$-spaces. In particular, a big problem  is the absence of a Wiener-Hopf factorization, which plays outstanding role in the study of Toeplitz plus Hankel operators on classical $H^p$-spaces.

\subsection{Spaces and operators.}
Let $l^p(\sZ)$ denote the complex Banach space of all sequences $\xi=(\xi_n)_{n\in\sN}$ of complex numbers with the norm $||\xi||_p= \left (\sum_{n\in\sZ} |\xi_n|^p \right )^{1/p}, \quad 1\leq p <\infty$. As usual, $\sZ$ denotes the set of all integers. If we replace $\sZ$ by the set of all non-negative integers $\sZ^+$, we get another Banach space $l^p(\sZ^+)$. It can be viewed as a subspace of $l^p(\sZ)$ and we will often write $l^p$ for $l^p(\sZ^+)$. By $P$, we now denote the canonical projection from $l^p(\sZ)$ onto $l^p(\sZ^+)$ and let $Q:=I-P$. Further, let $J$ refer to the operator on $l^p(\sZ)$ defined by
\begin{equation*}
 (J\xi)_n=\xi_{-n-1}, \quad n\in\sZ.
\end{equation*}
The operators $J,P$ and $Q$ are connected with each other by the relations
\begin{equation*}
J^2=I,\quad JPJ=Q,\quad JQJ=P.
\end{equation*}
For each $a\in L^p=L^p(\sT)$, let $(\widehat{a}_k)_{k\in\sZ}$  be the sequence of its  Fourier coefficients. The Laurent operator $L(a)$ associated with $a\in L^\infty(\sT)$ acts on the space $l^0(\sZ)$ of all finitely supported sequences on $\sZ$ by
\begin{equation}\label{Lau}
(L(a)x)_k:=\sum_{m\in\sZ} \widehat{a}_{k-m}x_m,
\end{equation}
where the sum in the right-hand side of \eqref{Lau} contains only a finite number of non-zero terms  for every $k\in\sZ$. We say that $a$ is a multiplier on $l^p(\sZ)$ if $L(a)x\in l^p(\sZ)$ for any $x\in l_0(\sZ)$ and if
 \begin{equation*}
||L(a)||:=\sup \{ ||L(a)x||_p: x\in l_0(\sZ), \, ||x||_p=1 \}
 \end{equation*}
is finite. In this case, $L(a)$ extends to a bounded linear operator on $l^p(\sZ)$, which is again denoted by $L(a)$. The set $M^p$ of all multipliers on $l^p(\sZ)$ is a Banach algebra under the norm $||a||_{M_p}:=||L(a)||$ --- cf.  \cite{BS:2006}. It is well-known that $M^2=L^\infty(\sT)$. Moreover, every function $a\in L^\infty(\sT)$ with bounded total variation $\mathrm{Var}\,(a)$ is in $M^p$ for every $p\in (1,\infty)$ and the Stechkin inequality
 \begin{equation*}
||a||_{M^p} \leq c_p(||a||_\infty+ \mathrm{Var}\,(a))
 \end{equation*}
holds with a constant $c_p$ independent of $a$. In particular, every trigonometric polynomial and every piecewise constant function on $\sT$ are multipliers on any space $l^p(\sZ)$, $p\in (1,\infty)$. By $C_p$ and $PC_p$ we, respectively, denote the closures of algebras of all trigonometric polynomials $\mathcal{E}$ and  all piecewise constant functions $PC$ in $M^p$. Note that $C_2$ is just the algebra $C(\sT)$ of all continuous functions on $\sT$, and $PC_2$ is the algebra $PC(\sT)$ of all piecewise continuous functions on $\sT$. We also note that the Wiener algebra $W$ of the functions with absolutely converging Fourier series is also a subalgebra of $M^p$ and
\begin{equation*}
W\subset C_p\subset PC_p\subset PC \quad\text{and}\quad M^p\subset L^\infty.
\end{equation*}
For this and other properties of multiplier cf.~\cite{BS:2006}. We also recall the equation $J L(a)J=L(\widetilde{a})$ often used in what follows.

Let $a\in M^p$. The operators $T(a): l^p\to l^p$, $x\mapsto PL(a)x$
and $H(a):l^p\to l^p$, $x\mapsto PL(a)Jx= PL(a)QJ$ are, respectively, called Toeplitz and Hankel operators with generating function $a$. It is well-known that $||T(a)||=||a||_{M^p}$ and $||H(a)||\leq||a||_{M^p}$ for every multiplier $a\in M^p$. In this section we also use the notation $T_p(a)$ or \/ $H_p(a)$ in order to underline that the corresponding Toeplitz or Hankel operator is considered on the space $l^p$ for a fixed $p\in(1,\infty)$. The action of the operators $T_p(a)$ and $H_p(a)$ on the elements from $l^p$ can be written as follows
\begin{align*}
T(a) &: (\xi_j)_{j\in \sZ_+} \to \left ( \sum_{k\in\sZ_+} \widehat{a}_{j-k}\xi_k  \right )_{j\in\sZ_+},   \\
H(a) &: (\xi_j)_{j\in \sZ_+} \to \left ( \sum_{k\in\sZ_+} \widehat{a}_{j+k+1}\xi_k  \right )_{j\in\sZ_+}.
\end{align*}

Let $GM^p$  denote the group of invertible elements in $M^p$.

\begin{lem}[cf. Ref.~\cite{BS:2006}]
  Let $p\in(1,\infty)$.
  \begin{enumerate}
  \item If $T_p(a)$ is Fredholm, then $a\in GM^p$.
    \item  If $a\in M^p$, then one of the kernels of the operators $T_p(a)$ or $T_q^*(a)$, $1/p+1/q=1$ is trivial.
    \item If $a\in GM^p$, the operator $T_p(a)$ is Fredholm, and $\mathrm{ind}\,T_p(a)=0$, then $T(a)$ is invertible on $l^p$.
  \end{enumerate}
\end{lem}

\subsection{Kernels of a class of Toeplitz plus Hankel\\ operators.}

The goal of this subsection is to present a method on how to study certain problems for Toeplitz plus Hankel operator $T_p(a)+H_p(b)$ defined on $l^p$ via known results obtained in Sections~\ref{sec2}-\ref{sec5}. Since $M^p\subset L^\infty(\sT)$, for any given elements $a,b\in M^p$ the operator $T_p(a)+H_p(b)\in \mathcal{L}(l^p)$ generates the operator $\mathbf{T}_s(a)+\mathbf{H}_s(b)\in \mathcal{L}(H^s)$, $1<s<\infty$ in an obvious manner. We denote by $\mathcal{L}^{p,q}(l^p)$ the collection of all Toeplitz plus Hankel operators acting on the space $l^p$ such that the following conditions hold:
\begin{enumerate}[label=(\alph*)]
  \item $a,b\in M^p$.
    \item $T_p(a)+H_p(b)\in \mathcal{L}(l^p)$ is Fredholm.
  \item $\mathbf{T}_q(a)+\mathbf{H}_q(b)\in \mathcal{L}(H^q)$, $1/p+1/q=1$ is Fredholm.
  \item Both operators acting on the spaces mentioned have the same Fredholm index.
\end{enumerate}
The following observation is crucial for what follows. The famous Hausdorff-Young Theorem connects the spaces $l^p$ and $H^q$, $1/p+1/q=1$ via Fourier transform $\mathcal{F}(a)=(\widehat{a}_n)_{n\in \sZ}$, $a\in H^q$.

\begin{thm}[Hausdorff \& Yang  {\cite{Ed:1982}}]\hfill{\phantom{asdfvgster}}
\begin{enumerate}[label=(\alph*)]\label{thm6.2}

  \item If $g\in H^p$ and $1\leq p\leq 2$, then $\mathcal{F} g\in l^q$ and
   $$
||\mathcal{F} g||_q\leq ||g||_p.
   $$
  \item If $\varphi=(\vp_n)_{n\in\sZ}\in l^p(\sZ)$ and  $1\leq p\leq 2$, then the series $\sum\vp_n e^{int}$ converges in $L^q$ to a function $\breve{\vp}$ and
      \begin{equation*}
||\breve{\vp}||_q\leq ||\vp||_p.
      \end{equation*}
\end{enumerate}
  \end{thm}

Theorem~\ref{thm6.2} gives rise to the following construction. Let $\widehat{H}^p$, $p\in (1,\infty)$ be the set of all sequences $(\chi_n)_{n\in \sZ_+}$ such that there exists a function $h\in H^p$ with the property $\mathcal{F} h=(\chi_n)_{n\in \sZ_+}$. The set $\widehat{H}^p$ equipped with the norm $|| (g_k)_{k\in \sZ_+}||:= ||g||_{H^p}$, becomes a Banach space isometrically isomorphic to $H^p$. Part (a) of the Hausdorff-Young Theorem assures that $\widehat{H}^p$ is densely continuously embedded in the space $l^q$ and part (b) claims that $l^p$, $1\leq p\leq 2$ is continuously and densely embedded into $\widehat{H}^q$.

\begin{thm}[VD \& BS \cite{DS:2018}]\label{thm6.3}
Let $p\in (1,\infty)$. If $T_p(a)+H_p(b)\in \mathcal{L}^{p,q}$, then the Fourier transform $\mathcal{F}$ is an isomorphism between  $\ker (\mathbf{T}_q(a)+\mathbf{H}_q(b))$ and $\ker (T_p(a)+H_p(b))$.
  \end{thm}
Let us give a sketch of the proof. The matrix representation $[\mathbf{T}_q(a)+\mathbf{H}_q(b)]$ of the operator $\mathbf{T}_q(a)+\mathbf{H}_q(b)$ in the standard basis $(t^n)_{n\in\sZ^+}$ induces a linear bounded operator on $\widehat{H}^q$, so that
 \begin{equation*}
[\mathbf{T}_q(a)+\mathbf{H}_q(b)]=(\widehat{a}_{i-j}+\widehat{b}_{i+j+1})_{i,j=0}^\infty
 \end{equation*}
 The operators $\mathbf{T}_q(a)+\mathbf{H}_q(b) \in \mathcal{L}(H^q)$ and $[\mathbf{T}_q(a)+\mathbf{H}_q(b)] \in \mathcal{L}(\widehat{H}^q)$ have the same Fredholm properties as $T_p(a)+H_p(b)\in \mathcal{L}(l^p)$ for $1<q\leq2$. Moreover, the operator $T_p(a)+H_p(b)\in \mathcal{L}(l^p)$ is the extension of $[\mathbf{T}_q(a)+\mathbf{H}_q(b)]$ with the same index. However, in this case 
 \begin{equation*}
\ker [\mathbf{T}_q(a)+\mathbf{H}_q(b)]=\ker (T_p(a)+H_p(b)),\quad T_p(a)+H_p(b) \in \mathcal{L}(l^p),
 \end{equation*}
and for $p\geq 2$ the assertion follows. The case $1<p\leq 2$ can be treated analogously.

Thus the main problem in using Theorem~\ref{thm6.2} is whether it is known that $T_p(a)+H_p(b)\in \mathcal{L}^{p,q}$. This is a non-trivial fact but the following result holds.
 \begin{prop}[VD \& BS \cite{DS:2018}]\label{prop6.4}
Let $a,b\in PC_p$, $1<p<\infty$. If the operator $T_p(a)+H_p(b) \in \mathcal{L}(l^p)$ is Fredholm, then $T_p(a)+H_p(b) \in\mathcal{L}^{p,q}$.
 \end{prop}
 The proof of this proposition can be carried out using ideas from \cite{RS:2012} and \cite{RSS:2011}. However, a simpler proof can be obtained if the generating function $a$ and $b$ satisfy the matching condition $a \widetilde{a}=b \widetilde{b}$ and the operators $T_p(c), T_p(d)$ are Fredholm say in $l^p$. Then $T_p(c), T_p(d)\in \mathcal{L}^{p,q}$ --- cf. \cite[Chapters 4 and 6]{BS:2006}, and using classical approach, which also works in $l^p$ situation, one obtains the result.

Thus, it is now clear how to extend the results of Sections~\ref{sec2}-\ref{sec5} to Toeplitz plus Hankel operators acting on spaces $l^p$. Let us formulate just one such result without going into much details.

 \begin{thm}\label{thm6.5}
Let $(a,b)\in PC_p\times PC_p$ be a Fredholm matching pair with the subordinated pair $(c,d)$, and let $\widehat{c}^{-1}_{+,j}$, $j\in \sZ_+$ be the Fourier coefficients of the function $c_+^{-1}$, where $c_+$ is the plus factor in the Wiener-Hopf factorization \eqref{WHF} of the function $c$ in $H^q$. If $\kappa_1:=\mathrm{ind}\, T_p(c)> 0$,  $\kappa_2:=\mathrm{ind}\, T_p(d)\leq 0$, then the kernel of the operator $T_p(a)+H_p(b)$ admits the following representation:
 \begin{enumerate}[label=(\alph*)]
      \item If $\kappa_1=1$ and $\sigma
(c)=1$, then
 $$\ker (T_p(a)+H_p(b))= \{0 \}$$

 \item If $\kappa_1=1$ and $\sigma
(c)=-1$, then
 $$\ker (T_p(a)+H_p(b))= \mathrm{lin}\, \mathrm{span}\{(\widehat{c}^{-1}_{+,j})_{j\in\sZ_+}
 \}.$$

   \item If $\kappa_1>1$ is odd, then
     \begin{align*}
 &\quad \ker (T_p(a)+H_p(b))\\
 &=\mathrm{lin}\, \mathrm{span}\{(\widehat{c}^{-1}_{+,j-(\kappa_1\!-\!1)/2-l}\!-\! \boldsymbol\sigma (c) \widehat{c}^{-1}_{+,j-(\kappa_1\!-\!1)/2+l})_{j\in\sZ_+}\colon l\!=\!0,\cdots, \frac{\kappa_1\!-\!1}{2} \}.
     \end{align*}

   \item If $\kappa_1$ is even, then
     \begin{align*}
&\quad\ker (T_p(a)+H_p(b)) \\
& =\mathrm{lin}\, \mathrm{span}\{(\widehat{c}^{-1}_{+,j-\kappa_1/2+l+1}\!-\! \boldsymbol\sigma
(c) \widehat{c}^{-1}_{+,j-\kappa_1/2-l})_{j\in\sZ_+}\colon l\!=\!0,1,\cdots, \kappa_1/2\!-\!1
 \}.
     \end{align*}
  \end{enumerate}
 \end{thm}

\begin{rem}
Sometime the study of invertibility of Toeplitz plus Hankel operators $T_p(a)+H_p(b)$ can be carried out even if it is not known, whether this operator belongs to $\mathcal{L}^{p,q}$. Thus $l^p$-versions of Theorem~\ref{thm2.9}, Proposition~\ref{prop3.9} and Corollary~\ref{cor3.10} can be directly proved.
 \end{rem}

\section{Wiener-Hopf plus Hankel Operators.\label{sec7}}

This section is devoted to Wiener-Hopf plus Hankel operators. Let  $\chi_E$ refer to the characteristic function of the subset $E\subset \sR$ and let $\mathcal{F}$ and $\mathcal{F}^{-1}$ be the direct and inverse Fourier transforms -- i.e.
\begin{equation*}
 \mathcal{F}\varphi(\xi):=\int\limits_{-\infty}^\infty e^{i\xi x} \varphi(x)\,dx,\quad
\mathcal{F}^{-1}\psi(x):=\frac1{2\pi}\int\limits_{-\infty}^\infty
e^{-i\xi x}\psi(\xi)\,d\xi,\;\;x\in\sR.
 \end{equation*}
In what follows, we identify the spaces $L^p(\sR^+)$ and
$L^p(\sR^-)$, $1\leq p \leq\infty$ with the subspaces $P(L^p(\sR))$ and $Q(L^p(\sR))$ of the space $L^p(\sR)$, where $P$ and $Q$ are the projections in $L^p(\sR)$  defined by $Pf:= \chi_{\sR^+} f$ and $Q:=I-P$, respectively.

Consider the set $G$ of functions $g$ having the form
\begin{equation}\label{Eq19}
  g(t)= (Fk)(t)+\sum_{j\in\sZ} g_j e^{i\delta_j t}, \quad t\in\sR,
 \end{equation}
where $k\in L(\sR)$,  $\delta_j\in \sR$, $g_j\in\sC$ and the series in the right-hand side of \eqref{Eq19} is absolutely convergent. Any function $g\in G$ generates an operator $W^0(a)\colon L^p(\sR)\to L^p(\sR)$ and
operators $W(g),H(g)\colon L^p(\sR^+)\to L^p(\sR^+)$ defined by
 \begin{equation*}
  W^0(g):=\mathcal{F}^{-1}g\mathcal{F}\varphi,  \quad
  W(g):=PW^0(g) , \quad
  H(g):= PW^0(g)QJ.
\end{equation*}
Here and throughout this section, $J\colon L^p(\sR)\to L^p(\sR)$ is the reflection operator defined by  $J\vp :=
\widetilde{\vp}$ and $\widetilde{\vp}(t):=\vp(-t)$ for any $\vp\in L^p(\sR)$, $p\in [1,\infty]$. Operators $W(g)$  and  $H(g)$ are, respectively, called Wiener-Hopf and Hankel operators on the semi-axis $\sR^+$.  It is well-known \cite{GF:1974} that for $g\in G$, all three operators above are bounded on any space $L^p$,  $p\in [1,\infty)$.

In particular, the operator $W(g)$ has the form
  \begin{equation*}
   W(g)\vp(t)  =\sum_{j=-\infty}^\infty g_j B_{\delta_j} \vp(t)
 +\int_{0}^\infty k(t-s) \vp(s)\,ds, \quad t\in\sR^+ ,
 \end{equation*}
where
 \begin{align*}
B_{\delta_j} \vp(t)&=\vp(t-\delta_j) \quad\text{if}\quad \delta_j\leq 0,\\
B_{\delta_j} \vp(t)&=\left \{
 \begin{array}{l}
 0, \quad 0\leq t \leq  \delta_j\\
 \vp(t-\delta_j), \quad t> \delta_j
 \end{array}
 \right.
  \quad \text{if}\quad \delta_j >0.
 \end{align*}
Moreover, for $g=\mathcal{F}k$ the operator $H(a)$ acts as
\begin{equation*}
 H(g)\vp(t)=\int_0^\infty k(t+s)\vp(s)\,ds
\end{equation*}
and for  $g=e^{\delta t}$, one has $H(a)\vp(t)=0$ if $\delta \leq0$ and
 \begin{align*}
 H(g)\vp(t) &=
     \left \{
   \begin{array}{ll}
    \vp(\delta-t),& \quad 0\leq t\leq \delta,\\
    0,& \quad t>\delta,
     \end{array}
    \right .
    \end{align*}
if $\delta >0$.

For various classes of generating functions, the Fredholm properties of operators $W(a)$ are well studied \cite{BS:2006,BKS:2002,CD:1969,Du:1973,Du:1977,Du:1979,GF:1974}. In particular, Fredholm and semi-Fredholm Wiener-Hopf operators are one-sided invertible, and for $a\in G$ there is efficient description of the kernels and cokernels of $W(a)$ and formulas for the corresponding one-sided inverses.

The study of Wiener-Hopf plus Hankel operators
 \begin{equation}\label{WHH}
 \sW(a,b)=W(a)+H(b), \quad a,b\in L^\infty(\sR).
 \end{equation}
is much more involved. The invertibility and Fredholmness  of such operators in the space $L^2$ is probably first considered by Lebre \emph{et al.} \cite{LMT:1992}. In particular, the invertibility of $\sW(a,b)$ has been connected with the invertibility of a block Wiener-Hopf operator $W(G)$. Assuming that $a$ admits canonical Wiener-Hopf factorization in $L^2$ and the Wiener-Hopf factorization of the matrix $G$ is known, Lebre \emph{et al.} provided  a formula for  $\sW^{-1}(a,b)$. Nevertheless, the difficulties with  Wiener-Hopf factorization of matrix $G$ influence the efficiency of the method.  For piecewise continuous generating functions, the conditions of Fredholmness are obtained in \cite{RSS:2011}. A different method, called equivalence after extension, has been applied to Wiener-Hopf plus Hankel operators $\sW(a,a)=W(a)+H(a)$ with generating function $a$ belonging to various functional classes \cite{CN:2009,CS:2011}. The Fredholmness, one-sided invertibility and invertibility of such operators  are equivalent to the corresponding properties of the Wiener-Hopf operator  $W (a \widetilde{a}^{-1})$. Therefore,  known results about the invertibility and Fredholmness of $W (a \widetilde{a}^{-1})$ can be retranslated to the operator $\sW(a,a)$. On the other hand, the equivalence after extension has not been used to establish any representations of the corresponding inverses.  Another approach has been exploited in \cite{KS:2000, KS:2001} to study  the essential spectrum and the index of the operators $I+H(b)$.

We now consider Wiener-Hopf plus Hankel operators \eqref{WHH} with generating functions $a,b\in G$ satisfying the matching condition \eqref{eqn1c}, where $\widetilde{a}(t)$  and $\widetilde{b}(t)$ denote the functions $a(-t)$ and $b(-t)$, respectively. Thus we now assume that
\begin{equation}\label{WHHmat}
a(t)a(-t)=b(t)b(-t)
\end{equation}
and define the subordinated functions $c$ and $d$ by
\begin{equation*}
\quad c(t):=a(t)b^{-1}(t), \quad d(t):=a(t)\widetilde{b}^{-1}(t)=a(t)b^{-1}(-t), \quad t\in \sR.
\end{equation*}
The assumption \eqref{WHHmat} allows us to employ the method developed in Sec\-ti\-ons~\ref{sec2}-\ref{sec5} and establish necessary and also sufficient conditions for invertibility and one-sided invertibility of the operators under consideration and obtain efficient representations for the corresponding inverses. Note that here we only provide the main results. For more details  the  reader can consult \cite{DS:2014b,DS:2017a,DS:2019b}.

Let  $G^+\subset G$ and $G^-\subset G$ denote the sets of functions, which  admit holomorphic extensions to the upper and lower
half-planes, respectively. If $g\in G$ is a matching function -- i.e. $g\widetilde{g}=1$, then according to \cite{DS:2014b,DS:2017a}, it can be represented in the form
\begin{equation*}
g(t)=\left(\boldsymbol\sigma(g)\,\widetilde{g}_+^{-1}(t)\right)e^{i\nu t}\Big(\frac{t-i}{t+i}\Big)^n g_+(t),
\end{equation*}
where $\nu=\nu(g)\in \sR$, $n=n(g)\in \sN$, $\boldsymbol\sigma(g)=(-1)^n g(0)$, $\widetilde{g}_+^{\pm1}(t)\in G^-$ and $\widetilde{g}_+(\infty)=1$. This representation is unique and the numbers $\nu(g)$ and $n(g)$ in the representation  \eqref{Eq19} are defined as follows:
 \begin{equation*}
    \nu(g):=\lim_{l\to\infty}\frac{1}{2l} [\arg g(t)]_{-l}^l,
\quad n(g):=\frac{1}{2\pi} [\arg
(1+g^{-1}(t)(\mathcal{F}k)(t)]_{t=-\infty}^\infty.
\end{equation*}
Moreover, following the considerations of Sections~\ref{sec2}-\ref{sec3}, for any right-invertible operator $W(g)$ generated by a matching function $g$, we can introduce  complementary projections $\mathbf{P}_g^{\pm}$ on $\ker W(g)$. More precisely, if $\nu<0$ or if $\nu=0$ and $n<0$, then $W(g)$ is right-invertible and the projections $\mathbf{P}_g^{\pm}$ are defined by
   \begin{equation*}
\mathbf{P}_g^{\pm}:=(1/2)(I\pm JQPW^0(g)P) \colon \ker W(g)\to \ker W(g).
   \end{equation*}
In addition, we also need the translation operator $\vp^+$ defined by the subordinated functions $c$ and $d$. Assume that $W(c)$ is right-invertible operator. Let $W_r^{-1}(c)$ be any of its right-inverses and consider the following function:
    \begin{align*}
\vp^+=\vp^+(a,b)&:=\frac{1}{2} (
W_r^{-1}(c)W(\tilde{a}^{-1})- JQW^0(c)P
W_r^{-1}(c)W(\tilde{a}^{-1}))\nn\\
&\quad + \frac{1}{2} JQW^0(\tilde{a}^{-1}),
  \end{align*}
We are ready to discuss the invertibility of Wiener-Hopf plus Hankel operators starting with necessary conditions in the case where at least one of the indices $\nu_1:=\nu(c)$ or $\nu_2:=\nu(d)$ is not equal to zero. The situation $\nu_1=\nu_2=0$ will be considered later on. Let $n_1$ and $n_2$ denote the indices $n(c)$ and $n(d)$, respectively.
 \begin{thm}[VD \& BS \cite{DS:2019b}]\label{thm7.1}
If $a,b\in G$, the operator $W(a)+H(b)$ is one-sided invertible in $L^p(\sR^+)$ and at least one of the indices $\nu_1$ or $\nu_2$ is not equal to zero, then:
  \begin{enumerate}[label=(\roman*)]
    \item  Either $\nu_1 \nu_2\geq 0$ or $\nu_1>0$ and $\nu_2<0$.
    \item If $\nu_1=0$ and $\nu_2>0$, then  $n_1 >-1$ or $n_1=-1$ and $\boldsymbol\sigma(c)=-1$.
    \item  If $\nu_1<0$ and $\nu_2=0$, then $n_2<1$ or $n_2=1$ and $\boldsymbol\sigma(d)=-1$.
  \end{enumerate}
  \end{thm}
  Consider now the case $\nu_1>0$ and $\nu_2<0$ in more detail. It can be specified as follows.
   \begin{thm}[VD \& BS \cite{DS:2019b}]\label{thmt.2}
 Let $\nu_1>0$, $\nu_2<0$, $n_1=n_2=0$ and $\fN_\nu^p$, $\nu>0$ denote the set of functions $f\in L^p(\sR^+)$ such that $f(t)=0$ for $t\in (0,\nu)$.

 \begin{enumerate}[label=(\roman*)]
   \item If the operator $W(a)+H(b)\colon L^p(\sR^+)\to L^p(\sR^+)$, $1<p<\infty$ is invertible from the left, then
     \begin{equation*}
  \vp^+(\mathbf{P}^+_d)\cap \fN_{\nu_1/2}^p= \{0\},
     \end{equation*}
  where $\vp^+=\vp^+(ae^{-i\nu_1 t/2}, be^{i\nu_1 t/2})$.

   \item If the operator $W(a)+H(b)\colon L^p(\sR^+)\to L^p(\sR^+)$, $1<p<\infty$ is invertible from the right, then
     \begin{equation}\label{rightinv}
  \vp^+(\mathbf{P}^+_{\overline{c}})\cap \fN_{-\nu_2/2}^p= \{0\},
     \end{equation}
where $\vp^+=\vp^+(\overline{a}e^{i\nu_2 t/2},\overline{\tilde{b}}e^{-i\nu_2 t/2})$.
 \end{enumerate}
   \end{thm}
Passing to the case $\nu_1=\nu_2=0$, we note that now the indices $n_1$ and $n_2$ take over.

  \begin{thm}[VD \& BS \cite{DS:2019b}]\label{thm7.3}

 Let $a,b\in G$, $\nu_1=\nu_2=0$  and the operator $W(a)+H(b)$ is invertible from the left. Then:

  \begin{enumerate}[label=(\roman*)]

 \item In the case $n_2\geq n_1$,  the index $n_1$ satisfies the inequality
     \begin{equation*}
     n_1\geq -1
     \end{equation*}
     and if $n_1=-1$, then $\boldsymbol\sigma(c)=-1$ and $n_2>n_1$.

  \item  In the case $ n_1> n_2$, the index $n_1$ satisfies the inequality
\begin{equation*}%\label{nec}
 n_1\geq 1,
 \end{equation*}
and the index $n_2$ is either non-negative or $n_2<0$ and $n_1\geq -n_2$.

 \end{enumerate}
  \end{thm}
The necessary conditions of the right-invertibility have similar form and we refer the reader to \cite{DS:2019b} for details. The proof of  the above results is based on the analysis of the kernel and cokernel of the operator $W(a)+T(b)$. In particular, if $W(a)+T(b)$ is left-invertible and one  of the corresponding  conditions is not satisfied, then this operator should have a non-zero kernel, which cannot be true.

The sufficient conditions of one-sided invertibility can be also formulated in terms of indices $\nu_1, \nu_2, n_1$ and $n_2$. For example, the following theorem holds.
\begin{thm}[VD \& BS \cite{DS:2019b}] Let $a,b\in G$ and indices $\nu_1, \nu_2, n_1$ and $n_2$ satisfy any of the following conditions:
  \begin{enumerate}[label=(\roman*)]
    \item $\nu_1<0$ and $\nu_2<0$.
    \item $\nu_1>0$, $\nu_2<0$, $n_1=n_2=0$, operator $W(a)+H(b)$ is normally solvable and satisfies the condition \eqref{rightinv}.
        \item $\nu_1<0$, $\nu_2=0$ and $n_2<1$ or $n_2=1$ and $\boldsymbol\sigma(d)=-1$.
    \item $\nu_1=0$, $n_1\leq0$ and $\nu_2<0$.
    \item $\nu_1=0$ and $\nu_2=0$
        \begin{enumerate}[label=(\alph*)]
          \item $n_1\leq0$, $n_2<1$;
          \item $n_1\leq 0$,  $n_2=1$ and $\boldsymbol\sigma(d)=-1$;
                        \end{enumerate}
  \end{enumerate}
 Then the operator $W(a)+H(b)$ is right-invertible.
  \end{thm}
The sufficient conditions for left-invertibility of $W(a)+H(b)$ can be obtained by passing to the adjoint operator. Here we are not going to discuss this problem in whole generality. However, we use assumptions, which allow to derive simple formulas for left- or right-inverses.  These conditions are not necessary for one-sided invertibility and the corresponding inverses can be also constructed even if the conditions above are not satisfied.

\begin{thm}[VD \& BS \cite{DS:2019b}]\label{thm7.5}
 Let $(a,b)\in  G \times G$ be a matching pair. Then:
  \begin{enumerate}
 \item If  $W(c)$ and $W(d)$ are left-invertible, the
operator $W(a)+H(b)$ is also left-invertible and one of its left-inverses has the form
  \begin{equation}\label{EqLI}
  (W(a)+H(b))_l^{-1}= W_l^{-1}(c) W(\tilde{a}^{-1})W_l^{-1}(d)(I - H(\tilde{d})) + W_l^{-1}(c)H(\tilde{a}^{-1}).
 \end{equation}
 \item If  $W(c)$ and $W(d)$ are right-invertible, the
operator $W(a)+H(b)$ is also right-invertible and one of its right-inverses has the form
  \begin{equation}\label{EqnRI}
  (W(a)+H(b))_r^{-1}\!=\! (I - H(\tilde{c}))W_r^{-1}(c) W(\tilde{a}^{-1})W_r^{-1}(d) + H(a^{-1})W_r^{-1}(d).\!\!
 \end{equation}
 \end{enumerate}
 \end{thm}
For  $g\in G$, the corresponding one-sided inverse of the operator $W(g)$ can be written by using Wiener-Hopf factorization of $g$ \cite{GF:1974}. We also note that if $W(c)$ and $W(d)$ are invertible, formulas \eqref{EqLI} and \eqref{EqnRI} can be used to write the inverse operator for $W(a)+H(b)$. They also play an important role when establishing inverse operators in a variety of situations not covered by Theorem~\ref{thm7.5}. The corresponding proofs run similar to considerations of Sections~\ref{sec4} and \ref{sec5}, but there are essential technical differences because the corresponding kernel spaces can be infinite dimensional.

It is worth noting that the generalized invertibility of  operators $W(a)+H(b)$ can also be studied --- cf. \cite{DS:2019b}.

\section{Generalized Toeplitz plus Hankel Operators}
Here we briefly discuss generalized Toeplitz plus Hankel operators. These operators are similar to classical Toeplitz plus Hankel operators considered in Section~\ref{sec2}, but the flip operator $J$ is replaced by another operator $J_\alpha$ generated by a linear fractional shift $\alpha$. It turns out that the classical approach of Sections~\ref{sec2}-\ref{sec5} can also be used but the application of the method is not straightforward and requires solving various specific problems. Therefore, in this section we mainly focus on the description of the kernels and cokernels of the corresponding operators. These results lay down foundation for the invertibility study. We also feel that the Basor-Ehrhardt method can be realized in this situation, but we are not going to pursue this problem here.

 The following construction is based on the considerations of \cite{DS:2016b}. Let $\sS$ denote the Riemann sphere. We consider a mapping $\alpha\colon \sS\to \sS$ defined by
 \begin{equation}\label{Rsph}
\alpha(z):=\frac{z-\beta}{\overline{\beta}z-1},
 \end{equation}
where $\beta$ is a complex number such that $|\beta|>1$.

Let us recall basic properties of $\alpha$.

 \begin{enumerate}
\item The mapping $\alpha\colon\sS^2\to \sS^2$ is  one-to-one,
$\alpha(\sT)=\sT$,  and if $\sD:=\{ z\in\sC:|z|<1\}$ is the interior
of the unite circle $\sT$ and $\overline{\sD}:=\sD\cup \sT$, then
\begin{equation*}%\label{eq5}
\alpha(\sD)=\sS^2\setminus \overline{\sD}, \quad
\alpha(\sS^2\setminus \overline{\sD})=\sD .
\end{equation*}
We note that $\alpha$ is an automorphism of the Riemann sphere
and the mappings $H^p \to \overline{H^p}$, $h\mapsto h\circ \alpha$
and $\overline{H^p} \to H^p$, $h\mapsto h\circ \alpha$ are
well-defined isomorphisms.

\item  The mapping $\alpha\colon\sT\to \sT$ changes the orientation
of $\sT$, satisfies the Carleman condition $\alpha(\alpha(t))=t$ for
all $t\in\sT$, and  possesses two fixed points -- viz.
 \begin{equation}\label{FP}
 t_+=(1+
\lambda)/\overline{\beta} \quad \text{and}\quad t_-=(1- \lambda)/\overline{\beta},
 \end{equation}
where $\lambda:=i \sqrt{|\beta|^2-1}$.

\item The mapping $\alpha$ admits the factorization
 \begin{equation*}
\alpha(t)=\alpha_+(t) \, t^{-1} \, \alpha_-(t)
\end{equation*}
with the factorization factors
\begin{equation*}
\alpha_+(t)=\frac{t-\beta}{\lambda}, \quad \alpha_-(t)=
\frac{\lambda t}{\overline{\beta} t-1}.
\end{equation*}

\item On the space $L^p$, $1<p<\infty$, the mapping $\alpha$ generates
a bounded linear operator $J_\alpha$, called weighted shift operator
and defined by
  \begin{equation*}
J_\alpha \vp(t):= t^{-1} \alpha_-(t) \vp(\alpha(t)), \quad t\in
 \sT.
 \end{equation*}
\end{enumerate}
Further, for $a\in L^\infty$ let $a_\alpha$ denote the composition
of the functions $a$ and $\alpha$ --- i.e.
 $$
 a_\alpha(t):=a(\alpha(t)), \quad t\in\sT.
 $$
The operators $J_{\alpha}$, $P$, $Q$ and $aI$ are connected with each other by the equations
  \begin{equation*}
 J_\alpha^2=I, \quad  J_\alpha a I=a_\alpha J_\alpha,\quad  J_\alpha P=Q J_\alpha, \quad  J_\alpha Q=P J_\alpha ,
 \end{equation*}
and for any $n\in \sZ$, one has $(a^n)_{\alpha}=(a_{\alpha})^n:
=a_{\alpha}^n$.
In addition to the Toeplitz operator $T(a)$, any element $a\in L^\infty$ defines another operator $H_\alpha:=P\alpha QJ_\alpha$, called generalized Hankel operator. Generalized Hankel operators are similar to classical Hankel operators $H(a)$. For example, analogously to \eqref{Widom} operators $H_\alpha$ are connected with Toeplitz operators by the relations
 \begin{align*}
   T(ab)=T(a)T(b)+H_\alpha(a)H_\alpha (b_\alpha), \\
    H_\alpha(ab)=T(a)H_\alpha(b)+H_\alpha(a)T(b_\alpha).
\end{align*}
On the space $L^p$, we now consider the operators of the form
$T(a)+H_\alpha(b)$ and call them generalized Toeplitz plus Hankel operators generated by the functions $a$, $b$ and the shift $\alpha$.

The classical approach of Section~\ref{sec2} can be also  employed to describe the kernels and cokernels of $T(a)+H_\alpha(b)$. To this aim we develop a suitable framework, which is not a straightforward extension of the methods of Section~\ref{sec2}. Let us now assume that $a$ belongs to the group of invertible elements $GL^\infty$ and the duo $a,b$ satisfy the condition
\begin{equation}\label{GenM}
a a_\alpha=b b_\alpha.
\end{equation}
Relation \eqref{GenM} is again called the matching condition and the duo $(a,b)$ are  $\alpha$-matching functions. To each matching pair, one can assign another $\alpha$-matching pair $(c,d):=(ab^{-1},ab_\alpha^{-1})$ called the subordinated pair for $(a,b)$. It is easily seen that $c c_\alpha=d d_\alpha=1$. In what follows , any element $g\in L^\infty$ satisfying the relation $gg_\alpha=1$ is referred to as $\alpha$-matching function. The functions $c$ and $d$
can also be expressed in the form $c=b_{\alpha}a_{\alpha}^{-1}$, $d= b_{\alpha}^{-1} a$. Besides, if $(c,d)$ is the  subordinated  pair for an
$\alpha$-matching pair $(a,b)$, then $(\overline{d},\overline{c})$
is the subordinated pair for the matching pair $(\overline{a},
\overline{b}_\alpha)$, which defines the adjoint operator
  \begin{equation*}
(T(a)+H_{\alpha}(b))^*=T(\overline{a})+H_{\alpha}(\overline{b}_\alpha).
 \end{equation*}
Rewrite  the operator $J_{\alpha}: L^p\mapsto L^p$ in the form
 $$
J_{\alpha} \vp(t):= \chi^{-1}(t) \vp (\alpha(t)),
 $$
where $\chi(t) = t/(\alpha_-(t))$, and note that if $\alpha$ is the shift \eqref{Rsph}, then:
 \begin{enumerate}
 \item $\chi\in H^\infty$ is an $\alpha$-matching function and $\wind \chi =1$, where $\wind \chi$ denotes the winding number of the function $\chi$ with respect to the origin.

 \item The function $\chi_\alpha \in \overline{H^\infty}$ and
$\chi_\alpha(\infty)=0$.

\item If $a,b\in L^\infty$ and $n$ is a positive integer, then
   \begin{equation*}
 T(a)+H_{\alpha}(b)=(T(a \chi^{-n})+H_{\alpha}(b \chi^n)) T(\chi^n).
 \end{equation*}
 \end{enumerate}
 These properties allow us to establish the following version of the Coburn-Simonenko Theorem for generalized Toeplitz plus Hankel operators.
 \begin{thm}\label{trm8.1}
Let $a\in GL^\infty$ and let $A$ denote any of the operators
$T(a)-H_{\alpha}(a \chi^{-1})$, $T(a)+H_{\alpha}(a \chi )$,
$T(a)+H_{\alpha}(a)$, $T(a)-H_{\alpha}(a)$. Then
\begin{equation*}
\min\{\dim \ker A, \dim\coker A\}=0.
\end{equation*}
 \end{thm}
 Further development runs similar to the one presented in Subsection~\ref{ssec2.3} and all results are valid if the operator $H(b)$ is replaced by $H_\alpha(b)$, $\widetilde{a}$ by $a_\alpha$ and $b$ by $b_\alpha$. For example, if $(a,b)$ is an $\alpha$-matching pair with the subordinated pair $(c,d)$, then the block Toeplitz operator $T(V(a,b))$ with the matrix function
 \begin{equation*}
V(a,b)=
\left(
  \begin{array}{cc}
    0 & d \\
    -c & a_\alpha^{-1} \\
  \end{array}
\right),
 \end{equation*}
can be represented in the form
 \begin{align*}
T(V(a,b))& =\left(%
\begin{array}{cc}
 0  & T(d) \\
  -T(c)  & T(a_{\alpha}^{-1}) \\
   \end{array}%
\right)   \nn\\[1ex]
 &= \left(%
\begin{array}{cc}
 -T(d)  & 0 \\
 0  & I \\
   \end{array}%
\right)
\left(%
\begin{array}{cc}
 0  & -I \\
  I  & T(a_{\alpha}^{-1}) \\
   \end{array}%
\right)
\left(%
\begin{array}{cc}
 -T(c)  &  0\\
  0  & I \\
   \end{array}%
\right).
\end{align*}
Moreover, assuming that  $T(c)$ is invertible from the right and $T_r^{-1}(c)$ is one of the right inverses, we can again represent the kernel of $T(V(a,b))$ as the direct sum
\begin{equation*}
\ker T(V(a,b))=\Omega(c) \dotplus \widehat{\Omega}(d),
\end{equation*}
 where
 \begin{align*}%\label{eq}
\Omega(c) &:=\left \{   (\vp,0)^T:\vp\in \ker T(c)\right\},\\
\widehat{\Omega}(d) &:=\left \{ (T_r^{-1}(c)
T(a_{\alpha}^{-1})s,s)^T:s\in \ker T(d) \right \}.
\end{align*}
Further, we observe that the operators
\begin{align*}
  \mathbf{P}_{\alpha}^{\pm}(c)&:=\frac{1}{2}\left ( I \pm J_\alpha QcP\right )\colon \ker T(c)\to\ker T(c),  \\
  \mathbf{P}_{\alpha}^{\pm}(d)&:=\frac{1}{2}\left ( I \pm J_\alpha QdP\right )\colon \ker T(d)\to\ker T(d),
\end{align*}
are projections on the corresponding spaces and consider the translation operators
\begin{equation*}
\vp_\alpha^{\pm}\colon \mathbf{P}_{\alpha}^{\pm}(d)\to \ker(T(a)\pm H_\alpha(b)),
\end{equation*}
which are defined similar to \eqref{eqnphi}, but with $J$ and $\widetilde{a}^{-1}$ replaced by $J_\alpha$ and ${a}_\alpha^{-1}$, respectively.

In this situation, Proposition~\ref{prop3.6} reads as follows.
 \begin{prop}
Assume that $(a,b)\in L^\infty \times L^\infty$ is a Fredholm matching pair.
If the operator $T(c)$ is right-invertible, then
\begin{align*}%\label{eq}
   \ker(T(a)+H_{\alpha}(b))& = \im \mathbf{P}_{\alpha}^-(c) \dotplus   \vp_{\alpha}^{+}(\im \mathbf{P}_{\alpha}^+(d)), \\
\ker(T(a)-H_{\alpha}(b))& =  \im \mathbf{P}_{\alpha}^+(c)      \dotplus     \vp_{\alpha}^{-}(\im
\mathbf{P}_{\alpha}^-(d)).
\end{align*}
\end{prop}
Thus there is a one-to-one correspondence between the spaces $\ker T(V(a,b))$ and $\ker\diag (T(a)+H_{\alpha}(b),T(a)-H_{\alpha}(b))$ and in order to establish
the latter, we need an efficient description of the spaces $\im \mathbf{P}_{\alpha}^+(c)$ and $\im \mathbf{P}_{\alpha}^+(d)$.

 Let $g$ stand for one of the functions $c$ or $d$. Then $g$ is a matching function and $T(g)$ is a Fredholm Toeplitz  operator. We may assume that $\ind T(g)>0$, since only in this case at least one of the spaces $\im \mathbf{P}_{\alpha}^+(g)$ or $\im \mathbf{P}_{\alpha}^-(g)$ is non-trivial.
 The following factorization result is crucial for the description of $\im \mathbf{P}_{\alpha}^{\pm}(g)$.

 \begin{thm}
If $g\in L^\infty$ satisfies the matching condition $g
g_{\alpha}=1$ and  $\wind g=n$, $n\in \sZ$, then $g$ can be represented in the  form
\begin{equation*}
g=\xi  g_{+}\chi^{-n}(g_{+}^{-1})_{\alpha}\,,
\end{equation*}
where $g_+$ and $n$ occur in the Wiener-Hopf factorization
 \begin{equation*}
g= g_-t^{-n} g_+, \quad g_-(\infty)=1,
 \end{equation*}
of the function $g$, whereas $\xi \in \{-1,1\}$ and is defined by
  \begin{equation}\label{eq19}
  \xi= \left ( \frac{\lambda}{\overline{\beta}}\right )^n g_+^{-1}\left ( \frac{1}{\overline{\beta}}\right
).
 \end{equation}
  \end{thm}
   \begin{defn}
The number $\xi$ in \eqref{eq19} is called the
$\alpha$-factorization signature, or simply, $\alpha$-signature of
$g$ and is denoted by $\boldsymbol\sigma_{\!\alpha}(g)$.
 \end{defn}

The $\alpha$-signature is used to describe the kernels of the operators $T(a)+H_\alpha(b)$ and for certain classes of generating functions it can be defined with relative ease. For example, assuming that the operator $T(g)$ is Fredholm, $n:=\ind T(g)$ and $g$ is H\"older continuous at the fixed point $t_+$ or $t_-$ of~\eqref{FP}, one can show that $\boldsymbol\sigma_{\!\alpha}(g)=g(t_+)$  or $\boldsymbol\sigma_{\!\alpha}(g)=g(t_-)(-1)^n$. For piecewise continuous functions, the situation is more complicated but still can be handled --- cf.~\cite{DS:2016b}.

  \begin{thm}
Let $g\in L^\infty$ be an $\alpha$-matching function such that the
operator $T(g):H^p\to H^p$ is Fredholm and $n:=\ind T(g)>0$. If
$g=g_-t^{-n} g_+$, $g_-(\infty)=1$ is the corresponding Wiener--Hopf
factorization of $g$ in $H^p$, then  the following systems of
functions $\cB_{\alpha}^{\pm}(g)$ form bases in the spaces $\im
\mathbf{P}_{\alpha}^{\pm}(g)$:
  \begin{enumerate}
    \item If $n=2m$, $m\in \sN$, then
    \begin{equation*}
\cB_{\alpha}^{\pm}(g):= \{g_+^{-1}(\chi^{m-k-1}\pm
\boldsymbol\sigma_{\!\alpha} (g)\chi^{m+k}):k=0,1,\cdots, m-1
 \},
    \end{equation*}
     and $\dim\im \mathbf{P}_{\alpha}^{\pm}(g)=m$.

    \item If $n=2m+1$, $m\in \sZ_+$, then
    \begin{align*}
&\cB_{\alpha}^{\pm}(g):= \{g_+^{-1}(\chi^{m+k}\pm
\boldsymbol\sigma_{\!\alpha} (g)\chi^{m-k}):k=0,1,\cdots, m
 \}\setminus \{0\},\\
   &  \dim\im \mathbf{P}_{\alpha}^{\pm}(g)=m +(1\pm\boldsymbol\sigma_{\!\alpha}(g))/2
     \,,
    \end{align*}
 and the zero element belongs to one of the sets
$\cB_{\alpha}^{+}(g)$ or $\cB_{\alpha}^{-}(g)$ --- viz. for $k=0$
one of the terms $\chi^m(1\pm\boldsymbol\sigma_{\!\alpha}(g))$ is
equal to zero.
   \end{enumerate}
\end{thm}

The above considerations provide a powerful tool for the study of  invertibility of generalized Toeplitz plus Hankel operators and it should be clear  that the results obtained in Sections~\ref{sec4} and \ref{sec5} can be extended to this class of operators. However, additional studies may be needed in certain situations.  Nevertheless, let us  mention  one of such results here.

\begin{prop} Assume that $(a,b)\in L^\infty\times L^\infty$ is a Fredholm matching pair and the operators $T(c)$ and $T(d)$ are right-invertible. Then $T(a)+H_\alpha(b)$ and $T(a)-H_\alpha(b)$ are also right-invertible and corresponding right inverses are given by
  \begin{equation*}
  (T(a)\pm H_\alpha(b))_r^{-1}= (I \mp H_\alpha(c_\alpha))T_r^{-1}(c) T(a_\alpha^{-1})T_r^{-1}(d) \pm H_\alpha(a^{-1})T_r^{-1}(d).
 \end{equation*}
 \end{prop}

In conclusion of this section, we note the works \cite{KLR:2007,KLR:2009} where more general operators with the shift~\eqref{Rsph} are considered. However, the conditions imposed on coefficient functions are more restrictive and the results obtained are less complete.

 \section{Final Remarks.}

We considered various approaches to the study of invertibility of Toeplitz plus Hankel operators and their close relatives. Before concluding this survey, we would like to mention two more problems of special interest. The first one is the construction of Wiener-Hopf factorizations for multipliers acting on the spaces $l^p$ and also for those on $L^p(\sR)$-spaces. Some ideas on how to proceed with this problem in the $l^p$-context have been noted in \cite{DS:2018}.

Another problem of interest is the study of the kernels and cokernels of Wiener-Hopf plus Hankel operators acting on $L^p(\sR)$-spaces, $p\neq 2$ and generated by multipliers from sets more involved than the set $G$ considered in Section~\ref{sec7}, such as piecewise continuous multipliers, for example.

In should be clear that the list of open problems in the theory of Toeplitz plus Hankel operators is not limited to those mentioned in this work and it is up to the interested reader to single out a one for further consideration.

\end{document}